\documentclass[11pt]{amsart}
\usepackage{latexsym,amsmath,amsfonts,amscd,amssymb}

\newcommand{\seq}[3]{{#1}_{#2}, \ldots, {#1}_{#3}}

\newcommand{\Y}{\Sigma \times {{\Bbb S}}^1}
\newcommand{\Spz}{\text{Sp}\, (2g,{\Bbb Z})}
\newcommand{\spinc}{\text{Spin}^{\Bbb C}}
\newcommand{\p}{\phi^{SW}}
\newcommand{\ima}{\sqrt{-1}}
\newcommand{\Ceh}{{\Bbb Q}[\eta,\theta]}

\newcommand{\rim}{\mathcal{R}\text{\it im}}

\newcommand{\la}{\langle}
\newcommand{\ra}{\rangle}
\newcommand{\surj}{\twoheadrightarrow}
\newcommand{\inc}{\hookrightarrow}
\newcommand{\ar}{\rightarrow}
\newcommand{\bd}{\partial}
\newcommand{\x}{\times}
\newcommand{\ox}{\otimes}
\newcommand{\iso}{\cong}

\newcommand{\CP}{{\Bbb C \Bbb P}}
\newcommand{\point}{\text{pt}}
\newcommand{\im}{\text{im}}

\newcommand{\Map}{\text{Map}}

\newcommand{\Diff}{\text{Diff}}

\newcommand{\Jac}{\text{Jac}}

\newcommand{\Gr}{\text{Gr}}

\newcommand{\GCD}{\text{GCD }}

\newcommand{\ind}{\text{ind}}
\newcommand{\PD}{\text{P.D.}}

\newcommand{\coker}{\text{coker}}

\newcommand{\cA}{\mathcal{A}}

\newcommand{\cC}{\mathcal{C}}

\newcommand{\cG}{\mathcal{G}}
\newcommand{\cH}{\mathcal{H}}
\newcommand{\cM}{\mathcal{M}}
\newcommand{\cI}{\mathcal{I}}

\newcommand{\cK}{\mathcal{K}}

\newcommand{\cL}{\mathcal{L}}

\newcommand{\cR}{\mathcal{R}}

\newcommand{\cZ}{\mathcal{Z}}

\renewcommand{\AA}{{\Bbb A}}
\newcommand{\CC}{{\Bbb C}}

\newcommand{\PP}{{\Bbb P}}
\newcommand{\QQ}{{\Bbb Q}}
\newcommand{\RR}{{\Bbb R}}

\renewcommand{\SS}{{\Bbb S}}
\newcommand{\ZZ}{{\Bbb Z}}

\renewcommand{\a}{\alpha}
\renewcommand{\b}{\beta}

\newcommand{\g}{\gamma}
\newcommand{\e}{\varepsilon}
\newcommand{\f}{\epsilon}
\newcommand{\h}{\theta}

\newcommand{\s}{\sigma}

\newcommand{\om}{\omega}
\newcommand{\q}{\psi}

\newcommand{\G}{\Gamma}

\renewcommand{\S}{\Sigma}
\newcommand{\D}{\Delta}
\renewcommand{\L}{\Lambda}
\renewcommand{\P}{\Phi}

\newcommand{\frs}{{\frak s}}

\theoremstyle{plain}
\begingroup
\newtheorem{thm}{Theorem}[section]
\newtheorem{cor}[thm]{Corollary}
\newtheorem{lem}[thm]{Lemma}
\newtheorem{prop}[thm]{Proposition}

\endgroup

\theoremstyle{definition}
\newtheorem{defn}[thm]{Definition}

\theoremstyle{remark}
\newtheorem{rem}[thm]{Remark}

\newtheorem{criterium}[thm]{Criterium}

\title[Seiberg-Witten-Floer homology of a surface times a
circle]{Seiberg-Witten-Floer homology of a surface times a circle
for non-torsion spin$^{\Bbb C}$ structures}

\author{Vicente Mu\~noz and Bai-Ling Wang}
\thanks{Key words: $4$-manifolds, Seiberg-Witten invariants,
Seiberg-Witten-Floer homology. \\
Mathematics Subject Classification. Primary: 57R57. Secondary: 57N13.}
\date{April, 1999. Revised May, 2002.}

\begin{document}
\renewcommand{\theenumi}{\alph{enumi}}

\begin{abstract}
  We determine the Seiberg-Witten-Floer homology groups of the $3$-manifold
  $\Y$, where $\S$ is a surface of genus $g \geq 2$, together with its
  ring structure, for a $\spinc$ structure with non-vanishing first Chern
  class. We give applications to computing Seiberg-Witten
  invariants of $4$-manifolds which are connected sums along surfaces and
  also we reprove the higher type adjunction inequalities obtained by
  Oszv\'ath and Szab\'o.
\end{abstract}

\maketitle

\section{Introduction}
\label{sec:1}

In this paper we study the gluing theory for Seiberg-Witten
invariants of $4$-manifolds split along the $3$-manifold $Y=\Y$,
where $\S$ is a surface of genus $g \geq 2$. This produces
applications to the determination of the Seiberg-Witten invariants
of $4$-manifolds which are constructed as connected sums of other
$4$-manifolds along embedded surfaces, and to obtain constrains
for the Seiberg-Witten invariants of $4$-manifolds containing an
embedded surface of genus $g$ and non-negative self-intersection.
The seminal work in this direction is provided by~\cite{taubes}
leading to a proof of the generalized Thom conjecture. Analysis of
this kind on non-trivial circle bundles over surfaces appears
in~\cite{mrowka}~\cite{thom}.

Before stating the results, we set up some notation. Let $X$ be a
closed, connected, oriented smooth $4$-manifold with $b^+>0$ and a
fixed homology orientation (i.e.\ an orientation of $H^1(X;\RR)
\oplus H^{2+}(X;\RR)$). For a $\spinc$ structure $\frs$, the
Seiberg-Witten
invariant~\cite{witten}~\cite{salamon}~\cite{morgan}~\cite{taubes1}~\cite{wang2}
is a linear functional
  $$
   SW_{X,\frs}: \AA(X) \ar \ZZ,
  $$
where $\AA(X)=\text{Sym}^*H_0(X) \ox \L^* H_1(X)$, the free graded
algebra generated by the class of the point
$x \in H_0(X)$ and the $1$-cycles $\g\in
H_1(X)$ (we understand rational coefficients).
We grade $\AA(X)$ by declaring the degree of $x$ to be $2$ and the
degree of the elements in $H_1(X)$ to be $1$.
The invariants are constructed by endowing $X$ with a metric $g$ and
studying the moduli space $\cM_{X}(\frs)$
of solutions $(A,\P)$ modulo gauge to the Seiberg-Witten equations
\begin{equation}
  \label{eqn:sw}
   \left\{ \begin{array}{l}
     \rho((F_A - \ima\,\xi)^+) =(\P \P^*)_0 \\
     D_A \P =0, \\
   \end{array} \right.
\end{equation}
where $\P$ is a section of the positive spin bundle $W^+$ of
$\frs$, $A$ is a connection on the determinant line bundle $L=\det
W^{\pm}$, $D_A:\G (W^+) \to \G(W^-)$ is the Dirac operator twisted
with the connection $A$, $\rho$ denotes Clifford multiplication,
$(\P \P^*)_0$ is the trace-free part of $(\P\P^*)$ interpreted as
an endomorphism of $W^+$, $\xi$ is a (small) closed real two-form
introduced as a perturbation.

Note that the invariants are zero on elements whose degree is not $d(\frs)$
where
  $$
  d(\frs)={ c_1(\frs)^2 - (2\chi(X) +3 \s(X)) \over 4}
  $$
is the dimension of $\cM_X(\frs)$.
When $b^+>1$ the Seiberg-Witten invariants are independent of metrics
and perturbations. When $b^+=1$ the Seiberg-Witten invariants
depend on a chamber structure. Fix a component $\cK_0$ of the positive cone
$\cK(X)= \{x \in H^2(X;\RR)-\{0\}/ x^2 \geq 0 \}$.
Then we say that the Seiberg-Witten
invariants are computed in $\cK_0$ when the metric $g$ and perturbation $\xi$
satisfy $(c_1(\frs)+\frac{1}{2\pi}[\xi]) \cdot
\om_g < 0$, where $\om_g \in H^2(X;\ZZ)$ is a period point for the metric $g$
lying in $\cK_0$.

A basic class for $X$ is a $\spinc$ structure with non-zero Seiberg-Witten
invariant.
By analogy with the definitions of simple type in the context of Donaldson
invariants~\cite[introduction]{hff}, we give the following
\begin{defn}
  Let $X$ be a $4$-manifold with $b^+>1$. We say that
\begin{itemize}
\item $X$ is of {\em simple type\/} if $SW_{X,\frs}(z)= 0$ for any
  $z$ in the ideal generated by $x$ in $\AA(X)$, for any $\frs$.
\item $X$ is of {\em $H_1$-simple type\/} if $SW_{X,\frs}(z)= 0$
  for any $z$ in the ideal of $\AA(X)$ generated by $H_1(X)$, for any $\frs$.
\item $X$ is of {\em strong simple type\/} if it is both of simple type and
  of $H_1$-simple type, i.e.\ $SW_{X,\frs}(z)= 0$ whenever $\deg(z)>0$, for any $\frs$.
\end{itemize}
\end{defn}

Note that when $X$ has $b_1=0$ it is automatically of $H_1$-simple type.
There are manifolds not of $H_1$-simple type (for instance any manifold
which is a connected sum $X \# \, \SS^1\x\SS^3$, where $X$ has $b^+>1$,
see~\cite[proposition 2.2]{OS}), but it is an open question whether all
$4$-manifolds with $b^+>1$ are of simple type.

Now we are ready to state our main result. On the one hand, we
have applications to computing the Seiberg-Witten invariants of
connected sums along surfaces (see~\cite{genusg}). We prove the
following results in section~\ref{sec:7}.

\begin{thm}
\label{thm:main1}
  Let $\bar X_1$ and $\bar X_2$ be $4$-manifolds with embedded surfaces
  $\S \inc \bar X_i$, $i=1,2$, of the same genus $g\geq 2$, self-intersection
  zero and representing non-torsion homology classes, and
  let $X=\bar X_1 \#_{\S} \bar X_2$ be their connected sum along $\S$.
  Suppose that $\bar X_1$, $\bar X_2$ are both of strong simple type, and let
  $\frs_1$, $\frs_2$ be $\spinc$ structures on  $\bar X_1$, $\bar X_2$
  respectively, such that $c_1(\frs_1) \cdot \S=c_1(\frs_2) \cdot \S\neq 0$ and
  $d(\frs_1)=d(\frs_2)=0$. Let $\frs_o$ be a $\spinc$ structure on $X$
  obtained by gluing $\frs_1$ and $\frs_2$ and let $z\in \AA(X)$ such
  that $d(\frs_o)=\deg z$. Then
$$
  \sum_{h\in\rim} SW_{X,\frs_o+h} (z) =\left\{
  \begin{array}{ll}
  SW_{\bar X_1,\frs_1} (1) \cdot SW_{\bar X_2,\frs_2} (1) \qquad & z=1,\;
  |c_1(\frs) \cdot \S | = 2g-2 \\
   0 & z=1, \; |c_1(\frs) \cdot \S | < 2g-2 \\ 0 & \deg(z)>0
  \end{array} \right.
$$
  where $\rim\subset H^2(X;\ZZ)$ is the subspace generated by
  the rim tori (cf.\ ~\cite{rim}).
  If any of the manifolds has $b^+=1$, then
  the Seiberg-Witten invariants are computed
  for the component of the positive cone containing $\e \PD[\S]$, where
  $\e=1$ if $c_1(\frs_o) \cdot \S <0$ and $\e=-1$ if
  $c_1(\frs_o) \cdot \S> 0$.

  Moreover if the connected sum is admissible
  (see definition~\ref{def:admissible}),
  then there are no basic
  classes $\frs$ of $X$ such that $0<|c_1(\frs) \cdot \S |<2g-2$, and
  the basic classes of $X$ with $c_1(\frs) \cdot \S =\pm(2g-2)$
  are in bijection with pairs of basic classes
  $(\frs_1, \frs_2)$ of $X_1$ and $X_2$ respectively, such that
  $c_1(\frs_1) \cdot \S=c_1(\frs_2) \cdot \S =\pm (2g-2)$.
\end{thm}

This theorem is analogous to~\cite[corollary 13]{genusg}
and~\cite[corollary 15]{genusg} in the case of Donaldson
invariants. It is generalised to the following analogue
of~\cite[theorem 9.5]{hff},

\begin{thm}
\label{thm:main2}
  Let $\bar X_1$ and $\bar X_2$ be $4$-manifolds with embedded surfaces
  $\S \inc \bar X_i$, $i=1,2$, of the same genus $g\geq 2$, self-intersection
  zero and representing non-torsion homology classes, and
  let $X=\bar X_1 \#_{\S} \bar X_2$ be their connected sum along $\S$.
  Suppose that $\bar X_1$, $\bar X_2$ are of $H_1$-simple type, and let
  $\frs_1$, $\frs_2$ be $\spinc$ structures on $\bar X_1$, $\bar X_2$
  respectively, such that $c_1(\frs_1) \cdot \S=c_1(\frs_2) \cdot \S\neq 0$.
  Let $\frs_o$ be a $\spinc$ structure on $X$
  obtained by gluing $\frs_1$ and $\frs_2$ and let $z\in \AA(X)$ such
  that $d(\frs_o)=\deg z$. Then
$$
 \sum_{h\in \rim} SW_{X,\frs_o+h} (z) =0
$$
  if $\deg(z)>0$ or if $c_1(\frs) \cdot \S \not\equiv 2g-2 \pmod 4$.
  If $X$ has $b^+=1$ then the Seiberg-Witten invariants are computed
  for the component of the positive cone containing $\e \PD[\S]$, where
  $\e=1$ if $c_1(\frs_o) \cdot \S <0$ and $\e=-1$ if
  $c_1(\frs_o) \cdot \S> 0$.

  Moreover if the connected sum is admissible
  then the basic classes $\frs$ of $X$ with $c_1(\frs) \cdot \S\neq 0$
  all satisfy $c_1(\frs) \cdot \S \equiv 2g-2\pmod 4$.
\end{thm}

The restriction $c_1(\frs) \cdot \S\neq 0$ in theorems
\ref{thm:main1} and \ref{thm:main2} is due to the fact that the
gluing theory for Seiberg-Witten invariants used here only works
for non-torsion $\spinc$ structures.

On the other hand, our analysis also leads to a different proof of
the higher type adjunction inequalities first obtained by
Oszv\'ath and Szab\'o in~\cite{OS}. Our method of proof is more
transparent and parallels that of~\cite{adjuncti} for proving the
higher type adjunction inequalities for Donaldson invariants.
Section~\ref{sec:8} is devoted to this issue.

\begin{thm}{\bf (\cite[theorem 1.4]{OS})}
\label{thm:main4}
  Let $X$ be a $4$-manifold and let $\S \subset X$ be an embedded surface of
  genus $g \geq 2$ representing a non-torsion homology class with
  self-intersection $\S^2\geq 0$. Let $a\in \AA(X)$ and $b\in \AA(\S)$.
  If $X$ has $b^+>1$ and $\frs$ is a $\spinc$ structure with
  $SW_{X,\frs}(a\, b) \neq 0$ and $|c_1(\frs)\cdot \S|  + \S^2>0$ then we have
$$
  |c_1(\frs)\cdot \S|  + \S^2 + \deg(b) \leq  2g-2.
$$
  Furthermore, when $b^+=1$ then for each $\spinc$ structure $\frs$ on $X$
  with $- c_1(\frs)\cdot \S  + \S^2 > 0$,
  for which $SW_{X,\frs}(ab) \neq 0$, when calculated in the component of
  $\cK (X)$ containing $\PD [\S]$, we have an inequality
$$
  - c_1(\frs)\cdot \S  + \S^2+ \deg (b) \leq 2g-2.
$$
\end{thm}

\begin{thm}{\bf (\cite[theorem 1.1]{OS})}
\label{thm:main5}
  Let $X$ be a $4$-manifold of $H_1$-simple type (e.g.\ with $b_1=0$)
  and let $\S \subset X$ be an embedded surface of genus $g\geq 2$
  representing a non-torsion homology class with $\S^2\geq 0$.
  If $X$ has $b^+>1$ then for each basic class $\frs$ for $X$
  with $|c_1(\frs)\cdot \S|  + \S^2>0$ we have
$$
  |c_1(\frs)\cdot \S|  + \S^2 + 2 d(\frs) \leq  2g-2.
$$
  If $b^+=1$ then for each basic class $\frs$ for the Seiberg-Witten
  invariants of $X$ calculated in the component of
  $\cK (X)$ which contains $\PD [\S]$ such that
  $- c_1(\frs)\cdot \S  + \S^2 > 0$,
  we have an inequality
$$
  - c_1(\frs)\cdot \S  + \S^2+ 2d (\frs) \leq 2g-2.
$$
\hfill$\Box$
\end{thm}

This is a particular case of the more general result

\begin{thm}{\bf (\cite[theorem 1.3]{OS})}
\label{thm:main6}
  Let $X$ be a $4$-manifold with an embedded surface $\imath:\S \inc X$
  of genus $g\geq 2$ representing a non-torsion homology class with
  $\S^2\geq 0$. Let $l$ be an integer so that there is a symplectic basis
  $\{\g_i\}_{i=1}^{2g}$ of $H_1(\S)$ with $\g_i\cdot\g_{g+i}=1$,
  $1\leq i \leq g$, satisfying that $\imath_* (\g_i)=0$ in $H_1(X)$
  for $i=1,\ldots, l$. Let $a\in \AA(X)$ and
  $b \in\AA(\S)$ be an element of degree $\deg(b) \leq l+1$.
  If $X$ has $b^+>1$ and $\frs$ is a $\spinc$ structure such that
  $SW_{X,\frs}(a\,b) \neq 0$ and $|c_1(\frs)\cdot \S|  + \S^2>0$ then we have
$$
  |c_1(\frs)\cdot \S|  + \S^2 + 2 \deg(b) \leq  2g-2.
$$
  Furthermore, when $b^+=1$ then for each $\spinc$ structure $\frs$ on $X$
  with $- c_1(\frs)\cdot \S  + \S^2 > 0$,
  for which $SW_{X,\frs}(a\,b) \neq 0$, when calculated in the component of
  $\cK(X)$ containing $\PD [\S]$, we have
$$
  - c_1(\frs)\cdot \S  + \S^2+ 2\deg (b) \leq 2g-2.
$$
\end{thm}

Again the restriction $|c_1(\frs)\cdot \S|  + \S^2>0$ is due to
the fact that we only use the Seiberg-Witten-Floer theory for
non-torsion $\spinc$ structures (see section \ref{sec:2}).

In order to prove these results,  we determine completely the
structure of the Seiberg-Witten-Floer (co)homology of the three
manifold $Y=\Y$, where $\S$ is a closed surface of genus $g\geq
2$, for any $\spinc$ structure with non-zero first Chern class. We
use the Seiberg-Witten-Floer groups as defined in~\cite{wang},
since they satisfy a gluing theorem for Seiberg-Witten invariants.
We prove

\begin{thm}
\label{thm:main7}
  Let $\frs_Y$ be a $\spinc$ structure on $Y=\Y$ with $c_1(\frs_Y) \neq
  0$. If $c_1(\frs_Y) \neq
  2r \PD[\SS^1]$, with $-(g-1) \leq r \leq g-1$ then $HFSW^*(Y,\frs_Y)=0$.
  Let $\frs_r$ be the $\spinc$ structure on $Y$ with $c_1(\frs_r)=
  2r \PD[\SS^1]$, $-(g-1) \leq r \leq g-1$, $r \neq 0$, and set $d=g-1-|r|$.
  Then there is an isomorphism of vector spaces
  $HFSW^*(Y,\frs_r) \iso H^*(s^d\S)$, where
  $s^d\S$ is the $d$-th symmetric product of $\S$.
\hfill $\Box$
\end{thm}

Theorem~\ref{thm:main7} will follow from theorem~\ref{thm:key} and
proposition~\ref{prop:isom-no}. The Seiberg-Witten-Floer homology
of $Y=\Y$ has a natural ring structure coming from the cobordism
which is a pair of pants times $\S$ (cf.\ ~\cite{donaldson}). This
should be closely related (if not isomorphic) to the quantum
cohomology of the symmetric products of $\S$ (see~\cite{thaddeus}
for a partial computation of the latter), in the same way as the
instanton Floer cohomology of $\Y$ is isomorphic to the quantum
cohomology of the moduli space of rank $2$ odd degree stable
vector bundles over $\S$ (see~\cite{munoz2}). We fully compute the
ring structure of $HFSW^*(Y,\frs_r)$, which is a deformation of
the ring structure on $H^*(s^{g-1-|r|}\S)$.

\begin{thm}
\label{thm:main8}
  Let $-(g-1)\leq r \leq g-1$, $r\geq 0$, and set $d=g-1-|r|$. Then
  there is a presentation
  $$
  HFSW^*(\Y,\frs_r)=\bigoplus_{k=0}^{d} \L^k_0 \ox
  \frac{\Ceh}{(\tilde\cR_k^g ,\h \tilde\cR_{k+1}^g,\eta^{d+1},\h^{d+1})},
  $$
  where $\h\in \L^2H_1(\S)$ represents the intersection form in $H^1(\S)$
  and
  $$
  \L^k_0=\L^k_0 H_1(\S)=\ker \left(\h^{g-k+1} :\L^k H_1(\S) \to \L^{2g-k+2}
  H_1(\S)\right)
  $$
  is the primitive part. The polynomials are defined as
  $$
   \tilde\cR_k^g= \sum_{i=0}^{\a}
   \frac{{d-k-\a+1 \choose i}}{i!{g-k \choose i}}(-1)^i \eta^{\a-i}\h^i -
   \sum_{i=0}^{\a+|r|} \frac{{\a+|r|\choose i}}{i!{g-k \choose i}}
   \eta^{\a+|r|-i}\h^i,
  $$
  where $\a=[{d-k \over 2}]+1$, $0\leq k \leq d$, and $\tilde\cR_{d+1}^g=1$.
\hfill $\Box$
\end{thm}

This theorem will follow from theorem~\ref{thm:cRk}.

\section{Review of Seiberg-Witten-Floer homology for non-torsion
$\spinc$ structures}
\label{sec:2}

Let $Y$ be an oriented $3$-manifold with first Betti number $b_1>0$ and
a $\spinc$ structure $\frs_Y$ on $Y$ with $c_1(\frs_Y)$ non-torsion.
We are going to review the construction of the Seiberg-Witten-Floer
(co)homology groups $HFSW_*(Y,\frs_Y)$ of $Y$ from~\cite{wang}.

\subsection{Definition of $HFSW(Y,\frs_Y)$}
\label{subsec:2.1}

We endow $Y$ with a metric $g_Y$ and fix a base connection $A_0$
on the determinant line bundle $L_Y=\det W_Y$, where $W_Y$ stands
for the spin bundle. There is a Chern-Simons Seiberg-Witten
functional (taking values in $\RR/\ZZ$) defined on the
configuration space of gauge classes $[A,\q]$ of a connection $A$
on $L_Y$ and a section $\q$ of $W_Y$,
 $$
  \cC_{\eta}(A,\q)=-{1\over 2}\left(\int_Y (A-A_0)\wedge (F_A+F_{A_0}- 2*_Y \ima
  \eta) + \la\q, D_A \q  \ra dvol_Y \right),
 $$
where $D_A$ stands for the Dirac operator on $W_Y$ coupled with
$A$, $\eta$ is a (perturbative) real coexact one-form on $Y$. The
critical points of $\cC_{\eta}$ correspond to translation
invariant solutions to~\eqref{eqn:sw} on the tube $Z=Y \x\RR$ for
the perturbation $\xi=* \eta$. The equations we are led to are
\begin{equation}
  \label{eqn:sw3d}
  \left\{ \begin{array}{l}
    {*}F_A = q(\q) + \eta \\
    D_A \q =0,
  \end{array} \right.
\end{equation}
where $q(\q)$ is the standard quadratic form. When $c_1(\frs_Y)$
is non-torsion we can choose a generic  perturbation parameter
$\eta$ to get a finite collection of non-degenerate irreducible
(oriented) points~\cite{wang}. Denoted by ${\cM}_Y(\frs_Y, \eta)$
be the set of solutions to the perturbed Seiberg-Witten
equations~\eqref{eqn:sw3d} on $(Y, \frs_Y)$. There is also a
well-defined relative index on ${\cM}_Y(\frs_Y, \eta)$, taking
values in $\ZZ/N\ZZ$, where
 $$
 N= \GCD \{ \la c_1(\frs_Y), \sigma\ra \, | \, \sigma \in H_2(Y,\ZZ)\}.
 $$
The Seiberg-Witten-Floer chain complex $CFSW_*(Y)$ is the vector space
generated by the gauge classes of solutions to~\eqref{eqn:sw3d} with the
relative grading.

The downward flow equation of $\cC_{\eta}$ is the 4-dimensional
Seiberg-Witten equation on $(Y \times \RR, g_Y + dt^2)$ with the
pull-back $\spinc$ structure under the temporal gauge:
\begin{equation}
  \label{eqn:sw4d}
  \left\{ \begin{array}{l}
    \displaystyle{ \frac{\partial A(t)}{\partial t}} =
    -  {*}F_A + q(\q) + \eta \\
   \displaystyle{ \frac{\partial \q (t)}{\partial t}} = -D_A \q .
  \end{array} \right.
\end{equation}
In general, in order to achieve the transversality property for
the moduli space of the Seiberg-Witten solutions on $Y \times \RR$
to~\eqref{eqn:sw4d} without destroying the $\RR$-translation
action and obtaining a natural compactification, we have to choose
a suitable perturbation of $\cC_\eta$ supported away from the set
of critical points. See~\cite{wang} for a detailed discussion.
(Similar kind of perturbations were first constructed
in~\cite{froyshov}.)

Let $a$ and $b$ in ${\cM}_Y(\frs_Y, \eta)$ be two gauge classes of
irreducible solutions to~\eqref{eqn:sw3d}. For generic
perturbations as in~\cite{wang}, any connected component of  the
moduli space $\hat{\cM}(a,b)$, the gauge classes of solutions to
the  perturbed Seiberg-Witten equations on the tube $Y\x \RR$ with
limits $a$ and $b$ respectively, is smooth, orientable and admits
a free $\RR$-action. We shall denote by $\cM^D(a,b)$ the
components of dimension $D$ in the quotient space of
$\hat{\cM}(a,b)$ by the $\RR$-action. Note that $D \equiv \ind (a,
b)-1 \pmod N$. We define a boundary map
\begin{eqnarray*}
  \bd:  CFSW_i(Y) & \ar & CFSW_{i-1}(Y) \\
        a & \mapsto &  \hspace{-5mm}
  \sum_{b\in \cM(Y, \frs_Y) \atop \ind(a, b)=1 \pmod N}\hspace{-5mm}
  \# \cM^0(a, b) \; b.
\end{eqnarray*}
The compactifications of the moduli space of the trajectory flow
lines ensure that $\bd$ is well-defined and
 $\bd^2=0$~\cite{wang} so we
obtain the Seiberg-Witten-Floer cohomology, denoted by
$HFSW^*(Y,\frs_Y)$, which is $\ZZ_N$-graded abelian group.

Notice that the first Chern class of the non-torsion $\spinc$
structure $\frs_Y$ defines a homomorphism $c_1(\frs_Y): H^1(Y,
\ZZ) \to \ZZ$ by $c_1(\frs_Y) ([u]) = \la [u]\cup c_1(\frs_Y),
[Y]\ra$ for any element $[u] \in H^1(Y, \ZZ)$. For any subgroup
$K\subset \ker (c_1(\frs_Y))$, there is a subgroup $\cG^K_Y$ of
the full gauge transformation group $\cG_Y$, whose elements lie in
the connected components determined by $K$ (since the group of
connected components of $\cG_Y$ is $H^1(Y, \ZZ)$). Consider the
Chern-Simons-Dirac functional on the configuration space $\cA_Y$
modulo the gauge group $\cG^K_Y$. The critical point set, denoted
by $\cM_{Y, K}(\frs_Y, \eta)$, is a covering space
 $$
 \pi_K: \cM_{Y, K}(\frs_Y, \eta) \longrightarrow \cM_Y (\frs_Y,
 \eta),
 $$
whose fiber is an $H^1(Y, \ZZ)/K$-homogeneous space.

There is a variant of the Seiberg-Witten-Floer chain complex whose
generators are elements in $\cM_{Y, K}(\frs_Y, \eta)$ with
relative $\ZZ$-graded indices and boundary operator $\partial^K$
given by counting the gradient flow lines of the perturted
Chern-Simons-Dirac functional on $\cA_Y/\cG^K_Y$ connecting two
critical points with relative index 1. We denote the
Seiberg-Witten-Floer homology in this setting by $HFSW_{*, [K]}(Y,
\frs_Y)$.

The following theorem was established in~\cite{wang} regarding
various properties of these Seiberg-Witten-Floer homologies of $(Y, \frs_Y)$.

\begin{thm} \label{Theorem:wang} (Theorem 1.1 \cite{wang})
For any closed oriented 3-manifold $Y$ with $b_1(Y)>0$ and a
$\spinc$ structure $\frs_Y$ such that $c_1(\frs_Y)\neq 0$ in
$H^2(Y, \ZZ)/\text{torsion}$, let $N$ be the divisibility of
$c_1(\frs_Y)$ in $H^2(Y, \ZZ)/\text{torsion}$,  that is
 $$
 N  = \GCD\{ \la c_1(\frs_Y), \sigma\ra \, | \, \sigma \in
 H_2(Y,\ZZ)\}.
 $$
Then there exists a finitely generated Seiberg-Witten-Floer
complex whose homology $HFSW_*(Y, \frs_Y)$ satisfies the following
properties:
\begin{enumerate}
\item $HFSW_*(Y, \frs_Y)$ is a topological invariant of $(Y, \frs_Y)$ and
 it is a $\ZZ_{N}$-graded abelian group.
\item There is an action of
 $$
  \AA (Y) = Sym^*(H_0(Y, \ZZ))\otimes \Lambda^*\bigl(H_1(Y,
   \ZZ)/\text{torsion}\bigr)
 $$
 on $HFSW_*(Y, \frs_Y)$ with elements in $H_0(Y, \ZZ)$ and $H_1(Y,
 \ZZ)/\text{torsion}$ decreasing degree in $HFSW_*(Y, \frs_Y)$ by
 $2$ and $1$ respectively.
\item For $(-Y, -\frs_Y)$, where $-Y$ is $Y$ with the reversed
 orientation and $-\frs_Y$ is the induced $\spinc$ structure, the
 corresponding Seiberg-Witten-Floer complex $C_*(-Y, -\frs_Y)$ is
 the dual complex of $C_*(Y, \frs_Y)$. There is a natural pairing
 $$
  \la \ , \ \ra : HFSW_*(Y, \frs_Y) \times HFSW_{-*}(-Y, -\frs_Y)
  \longrightarrow \ZZ
 $$
 such that $<z\cdot \Xi_1, \Xi_2> = <\Xi_1, z\cdot\Xi_2>$ for any
 $z\in \AA (Y)\cong \AA(-Y)$ and any cycles $\Xi_1 \in HFSW_*(Y,
 \frs_Y)$ and $\Xi_2\in HFSW_{-*}(-Y, -\frs_Y)$ respectively.
\item For any subgroup $K\subset \ker (c_1(\frs_Y))\subset H^1(Y, \ZZ)$,
 there is a variant of Seiberg-Witten-Floer homology denoted by
 $HFSW_{*, [K]}(Y, \frs_Y)$, which is a topological invariant and a
 $\ZZ$-graded $\AA (Y)$ module. For any $[u] \in H^1(Y, \ZZ)/ K$,
 there is an action of $[u]$ on $HFSW_{*, [K]}(Y, \frs_Y)$
 decreasing degrees by $\la [u]\wedge c_1(\frs_Y), [Y]\ra$. There
 is natural pairing
\begin{equation*}
 \la \ , \ \ra : HFSW_{*, [K]}(Y, \frs_Y) \times
 HFSW_{-*, [K]}(-Y, -\frs_Y)\longrightarrow \ZZ
\end{equation*}
 satisfying $<z\cdot\Xi_1, \Xi_2> = <\Xi_1, z\cdot\Xi_2>$ for any
 $z\in \AA (Y)\cong \AA(-Y)$ and any cycles $\Xi_1 \in HFSW_{*,
 [K]}(Y, \frs_Y) $ and $\Xi_2 \in HFSW_{-*, [K]}(-Y, -\frs_Y)$
 respectively. There is a $\AA (Y)$-equivariant homomorphism:
  $$
  \pi_K: HFSW_{*, [K]}(Y, \frs_Y)  \longrightarrow HFSW_{*}(Y,
  \frs_Y).
  $$
 If $K_1\subset K_2$ are two subgroups in $\ker (c_1(\frs_Y))$,
 there is a $\AA (Y)$-equivariant homomorphism $HFSW_{*, [K_1]}(Y,
 \frs_Y)  \to HFSW_{*, [K_2]}(Y, \frs_Y)$. Moreover, for any $m\in \ZZ$,
  \begin{equation*} 
  \pi_{\ker(c_1(\frs_Y))}:
  HFSW_{m, [\ker(c_1(\frs_Y))]} (Y, \frs_Y)  \cong HFSW_{m \pmod N}
  (Y, \frs_Y).
  \end{equation*}
\end{enumerate}
\end{thm}

\subsection{Relative Seiberg-Witten invariants and gluing formula}
\label{subsec:2.3}

Let $X_1$ be an oriented, connected $4$-manifold furnished  with a
cylindrical end of the form $Y \x [0,\infty)$. Suppose we have a
$\spinc$ structure $\frs_1$ over $X_1$ whose restriction to $Y$ is
a non-torsion $\spinc$ structure $\frs_Y$. Consider finite energy
solutions to the Seiberg-Witten equations on $X_1$ with finite
variations of the perturbed Chern-Simons-Dirac functional on the
end, there is an associated boundary asymptotic value map
 $$
 \partial_\infty: \cM_{X_1}(\frs_1) \to \cM_{Y, X_1}(\frs_Y, \eta)
 $$
where $\cM_{Y, X_1}(\frs_Y, \eta)$ is the quotient of solutions to
the perturbed Seiberg-Witten equations on $(Y, \frs_Y)$ by the
action of those gauge transformations which can be extended to
$X_1$. In fact, $\pi_1 : \cM_{Y, X_1}(\frs_Y, \eta) \to
\cM_Y(\frs_Y, \eta)$ is a covering map with fiber an $H^1(Y,
\ZZ)/\im (i_1^*)$-homogeneous space. Here $\im (i_1^*) \subset
\ker (c_1(\frs_Y))$ is the image of the map $i^*_1 : H^1(X_1, \ZZ)
\to H^1(Y, \ZZ)$ induced from the boundary embedding map $i_1$.
Generically, the fiber of $\partial_\infty$ is an oriented, smooth
manifold of dimension given by Atiyah-Patodi-Singer index theorem,
and it can be compactified to a smooth manifold with corners.
See~\cite{wang} for the detailed discussion.

The relative Seiberg-Witten invariant of $(X_1, \frs_1)$, as
defined in~\cite{wang}, takes values in the Seiberg-Witten-Floer
homology $HFSW_{*, [\im(i_1^*) ]}(Y, \frs_Y)$, and defines an
$\AA(Y)$-equivariant linear map
 $$
 \p_{X_1}(\frs_1,\cdot):
 \AA(X_1) \longrightarrow HFSW_{*, [\im(i_1^*) ]}(Y, \frs_Y).
 $$
Here the $\AA(Y)$-action on $\AA(X_1)$ is induced from the
homomorphism $(i_1)_*: \AA(Y)\to \AA(X_1)$. For any $z_1 \in
\AA(X_1)$ of degree $d$, $\p_{X_1}(\frs_1, z_1)$ can be expressed
in terms of the Seiberg-Witten invariant from the components of
dimension $d$ in $\cM_{X_1}(\frs_1)$.

For a second $4$-manifold $X_2$ with a cylindrical end $(-Y) \x
[0,\infty)$, we construct $X=X_1 \cup_Y X_2$ by cutting the ends
and gluing along the common boundary $Y$. The resulting manifold
may depend on the isotopy class of the diffeomorphism identifying
the boundaries, but we shall not make the dependence explicit. If
there is a $\spinc$ structure $\frs_2$ on $X_2$ with $\frs_2|_Y
\iso \frs_Y$, then we can glue $\frs_1$ and $\frs_2$. The
indeterminacy for the gluing is parametrised by $\coker
(H^1(X_1;\ZZ) \oplus H^1(X_2;\ZZ) \ar H^1(Y;\ZZ))$. The following
gluing formula is taken from~\cite{wang}.

\begin{thm}[Theorem 1.2 \cite{wang}]
\label{thm:gluing}
 Let $X$ be a closed manifold with $b^+ \geq 1$ which is written as
 $X=X_1 \cup_Y X_2$, where $X_1$ and $X_2$ are $4$-manifolds with
 boundary and $\bd X_1=-\bd X_2=Y$. Suppose that we have $\spinc$
 structures $\frs_1$ and $\frs_2$ on $X_1$ and $X_2$ respectively
 such that $\frs_1|_Y \iso \frs_2|_Y \iso \frs_Y$ is a non-torsion
 $\spinc$ structure on $Y$. Then for any $\spinc$ structure
 $\frs =\frs_1\#_{[u]}\frs_2$ obtained by gluing $\frs_1$ and
 $\frs_2$ along $Y$ using an isomorphism $u\in \Map (Y, U(1))$
 representing $[u] \in \coker (H^1(X_1;\ZZ) \oplus H^1(X_2;\ZZ) \ar
 H^1(Y;\ZZ))$, we have the following gluing formula for $z_i \in
 \AA(X_i)$, $i=1,2$,
  $$
    SW_{X,\frs} (z_1z_2)= \la [u]( \pi_1(\p_{X_1}(\frs_1,z_1))),
    \pi_2(\p_{X_2}(\frs_2,z_2))\ra.
  $$
 Here  $[u]$ acts on $HFSW_{*, [\im(i^*_1)+\im(i^*_2)]}(Y,
 \frs_Y)$, $\pi_1$ and $\pi_2$ are the $\AA(Y)$-equivariant
 homomorphisms induced from the inclusion maps $\im(i^*_1)\subset
 (\im(i^*_1)+\im(i^*_2))$ and $\im(i^*_2)\subset (\im(i^*_1)+
 \im(i^*_2))$ respectively, and the pairing on the right hand side
 is the natural pairing
  $$
  HFSW_{*, [\im(i^*_1)+\im(i^*_2)]}(Y, \frs_Y)
   \times HFSW_{*, [\im(i^*_1)+\im(i^*_2)]}(-Y, -\frs_Y),
  $$
 with the degrees in $HFSW_{*, [\im(i^*_1)+\im(i^*_2)]}(-Y,
 -\frs_Y)$ shifted by $\deg(z_1)+\deg (z_2)$. When $b^+=1$, the
 Seiberg-Witten invariants correspond to a metric giving a long
 neck. In particular, let ${\mathcal S}$ be the set of $\spinc$
 structures on $X$ with $\frs|_{X_i}=\frs_i$, $i=1,2$ and $\frac14
 (c_1(\frs)^2 -2(\chi(X) + \sigma (X))) = \deg (z_1)+\deg (z_2)$,
 then
 $$
   \sum_{\frs \in {\mathcal S}} SW_{X,\frs} (z_1z_2)=
   \la \pi( \p_{X_1}(\frs_1,z_1)), \pi (\p_{X_2}(\frs_2,z_2))\ra.
  $$
  Here $\pi(\p_{X_i}(\frs_i, z_i))$ are elements in
  $HFSW_{*}(\pm Y, \pm \frs_Y)$ under the maps
  $$
  \pi: HFSW_{*, [\im(i^*_1)+\im(i^*_2)]}(\pm Y, \pm \frs_Y)
   \to HFSW_{*}(\pm Y, \pm \frs_Y)
   $$
  and the pairing on the right hand side is the pairing
  $HFSW_{*}(Y, \frs_Y) \times HFSW_{*}(-Y, -\frs_Y) \to \ZZ$.
\end{thm}

\section{Seiberg-Witten-Floer homology of $\Y$}
\label{sec:4}

  From now on we shall consider the three-manifold $Y=\Y$, which is
the central object of our study. As $H^2(Y;\ZZ)$ has no
$2$-torsion, the $\spinc$ structures $\frs_Y$ on $Y$ are
determined by the determinant line bundle $L_Y=c_1(\frs_Y)$. As
$c_1(\frs_Y)$ reduces to $w_2(Y)=0$ modulo $2$, it has to be an
even class in $H^2(Y;\ZZ)$.

\begin{prop}
\label{prop:sw3d-Y}
  Let $\frs_Y$ be a $\spinc$ structure on $Y$.
  Let $\cM$ be the moduli space of solutions to~\eqref{eqn:sw3d}
  with zero perturbation. Then $\cM$ is empty unless $c_1(\frs_Y)=
  2r \PD[\SS^1]$, with $-(g-1) \leq r \leq g-1$. For $c_1(\frs_Y)=
  2r \PD[\SS^1]$, with $-(g-1) \leq r \leq g-1$ and $r\neq 0$, $\cM$
  is Morse-Bott irreducible and isomorphic to $s^d \S$ with $d=g-1-|r|$.
\end{prop}

\noindent {\em Proof.}
  We choose a rotation invariant metric for $Y$ of the form
  $g_{\S} + d\h \otimes d\h$, where $g_{\S}$ is a metric on
  $\S$ with unit area and scalar curvature $-4\pi (2g -2)$,
  and $\h$ is the coordinate on $\SS^1=\RR/\ZZ$.
  Think of $\Y$ as $\S \times [0,1]$ with the boundaries $\S \times \{0\}$
  and $\S \times \{1\}$ identified by the identity. The line bundle $L_Y$
  is constructed out of the pull back under the projection $\S \times [0,1]
  \ar\S$ of a line bundle $L_{\S}$ on $\S$ by gluing along the boundaries
  with an isomorphism $\s \in \cG_{\S}= \Map (\S, \SS^1)$. Then $c_1(L_Y)=
  c_1(L_{\S}) + [\s] \ox [\SS^1]$, where $[\s]$ is the class of $\s$ in
  $[\S ; \SS^1 ] \iso H^1(\S;\ZZ)$. The $\spinc$ structure $\frs_Y$ induces
  a $\spinc$ structure $\frs_{\S}$ on $\S$ with determinant
  line bundle $L_{\S}$. The spin bundle of $\frs_{\S}$ is of the form
  $W_{\S} = (\L^0 \oplus \L^{0,1}) \otimes \mu$, for a line bundle
  $\mu$ such that $L_{\S}= K_{\S}^{-1} \otimes \mu^2$, where
  $K_{\S}$ stands for the canonical bundle of $\S$.

  Now consider any solution $(A,\q)$ to~\eqref{eqn:sw3d} with $\eta=0$.
  In $\S \x [0,1]$, we can kill, by a gauge transformation, the $d\h$
  component of $A$, i.e.\ we can suppose that we have a family $A_{\h}$,
  $\h \in [0,1]$, of connections on $L_{\S}$ (up to a constant gauge)
  with the boundary condition $A_1 = \s^*(A_0)$, for some $\s \in \cG_{\S}$
  in the homotopy class determined by $L_Y$. So $(A,\q)$ is interpreted as
  a path $(A_{\h},\a_{\h},\b_{\h})$,  $\h \in [0,1]$,
  where $\a_{\h} \in \L^0 \otimes \mu$ and $\b_{\h} \in \L^{0,1} \otimes \mu$.
  Let us rewrite equations~\eqref{eqn:sw3d} in this set-up. Clearly
  $*F_A = \L F_{A_{\h}} d\h + *_{\S}( {\bd A \over \bd \h})$,
  and the map $q$ and the Dirac operator are as follows
  $$
    q(\q)= -* {\left( \a \overline{\b}  + {\overline{\a}} \b\right) \over
   2 } + {|\a|^2-|\b|^2 \over 2} d\h,
  $$
  $$
  D_A= \left(\begin{array}{cc} -\ima {\bd \over \bd \h} & \sqrt{2}
 \bar \bd_A^* \\
  \sqrt{2}\bar \bd_A & -\ima {\bd \over \bd \h} \end{array} \right) :
  (\L^0\ox\mu) \oplus (\L^{0,1} \ox \mu) \ar
  (\L^0\ox\mu) \oplus (\L^{0,1} \ox \mu).
  $$
  So the solutions to~\eqref{eqn:sw3d} correspond to solutions to
  \begin{equation}
    \left\{ \begin{array}{l}
      {\bd \a \over \bd \h }= -\ima \sqrt{2} \bar \bd_{A_{\h}}^* \b \\
      {\bd \b \over \bd \h }= \ima \sqrt{2} \bar \bd_{A_{\h}} \a \\
      2{\bd A_{\h} \over \bd \h }= -\ima (\a \overline{\b} +\b \overline{\a}) \\
      2 \, \ima \L F_{A_{\h}} = -| \a |^2 + |\b|^2
    \end{array} \right. \label{eqn:sw2d}
  \end{equation}
  We can write $A_{\h}= \bd_{A_{\h}} + \bar \bd_{A_{\h}}$, so
  the third line is ${\bd \over \bd \h } (\bd_{A_{\h}})= -{\ima \over 2} \a
  \overline{\b}$, ${\bd \over \bd \h } (\bar \bd_{A_{\h}})= -{\ima \over 2}
  \overline{\a} \b$. Now suppose we have a solution to~\eqref{eqn:sw2d}. Then
  we work out the following expression (using $\ima \bar \bd^* = \L \bd$ on
  $(0,1)$-forms and $|\b|^2= - \ima \L \b \wedge \overline{\b}$)
$$
  {\bd \over \bd \h}(\bar \bd^* \b ) = -{\bd \over \bd \h}(\ima\L\bd) \b + \bar
  \bd^*  ({\bd\b \over \bd\h}),
 $$
  with the given equalities to get
 $$
  -{1 \over  \sqrt{2}\ima}{\bd ^2 \a \over \bd \h^2}= \ima\L{\ima \over
  2} \a \overline{\b}\b -\bar \bd^*( -\ima \sqrt{2} \bar\bd \a),
 $$
 $$
  -{\bd ^2 \a \over \bd \h^2} + {\sqrt{2} \over 2}
      \a |\b|^2 +2\bar\bd^*\bar\bd\a =0.
 $$
  Take scalar product with $\a$ and integrate along $\S$ by parts to get
 $$
  - \int_{\S} \la {\bd ^2 \a \over \bd \h^2} , \a \ra
  +  {\sqrt{2} \over 2}\int_{\S} | \a|^2 |\b|^2 +2 \int_{\S} |\bar \bd \a|^2=0,
 $$
  for every $\h \in [0,1]$. This equation makes sense in $\SS^1$, since
  the values for $\h=0$ and $\h=1$ coincide. Then we can integrate again by parts
  to get
 $$
  ||{\bd \over \bd \h} \a||^2 + {\sqrt{2} \over 2}
  || \a \b||^2+2 ||\bar \bd \a||^2 =0.
 $$
So either $\a=0$ or $\b=0$. In any case, $A_{\h}$, $\a_{\h}$
and $\b_{\h}$ are constant, i.e.\ if the line bundle $L_Y$ admits
solutions to~\eqref{eqn:sw3d} then it is pulled-back from $\S$ and
any solution is invariant under rotations in the $\SS^1$ factor.

Assume now that $c_1(L_Y)=2r \PD [\SS^1]$. For any solution
to~\eqref{eqn:sw2d}, either $\a =0$, $\bar \bd_{A_0}^* \b=0$ or
$\b =0$, $\bar  \bd_{A_{0}} \a=0$. Also $2r = c_1(L_{\S}) ={\ima\over
2\pi} \int_{\S} F_A={1\over 4\pi}
 \int_{\S} (|\b|^2-|\a|^2)$. If $r <0$ then $\b=0$ and
the solutions to equations~\eqref{eqn:sw2d} are equivalent to the solutions to
   $$
  \left\{ \begin{array}{l}
   \bar \bd_A \a =0 \\
    2\ima \L F_A = - |\a|^2
  \end{array} \right. \label{eqn:vortex}
   $$
on $\S$. These are the typical vortex equations. The space of
solutions is $s^d \S$, where $d=g-1+r$. If $r<-(g-1)$ then there
are no solutions. The case $r>0$ is analogous. \hfill $\Box$
\hspace{2mm}

\begin{thm}
\label{thm:key}
  Let $\frs_Y$ be a $\spinc$ structure on $Y$ with $c_1(\frs_Y)\neq 0$.
  Then $HFSW^*(Y,\frs_Y)=0$ unless $c_1(\frs_Y)=
  2r \PD[\SS^1]$, with $-(g-1) \leq r \leq g-1$.
  Let $\frs_r$ be the $\spinc$ structure on $Y$ with $c_1(\frs_r)=
  2r \PD[\SS^1]$,  $-(g-1) \leq r \leq g-1$, $r\neq 0$.
  Put $d=g-1-|r| \geq 0$, then $\dim HFSW^*(Y,\frs_r) \leq \dim H^*(s^d\S)$.
\end{thm}

\noindent {\em Proof.}
  The first claim is a direct consequence of proposition~\ref{prop:sw3d-Y}.
  Also from proposition~\ref{prop:sw3d-Y}, we know that
  the unperturbed Chern-Simons Seiberg-Witten functional already has
  non-degenerate critical manifolds. As in~\cite{Fukaya} and
  \cite[proposition 6]{braam-donaldson}, we can choose a perturbation
  modelled on the finite dimensional critical manifold $s^d\S$.
  Choose a positive and  perfect Morse function $f$ on $\S$, i.e.\ $f$ has one
  critical point of index 0, $2g$ critical points of index $1$ and one
  critical point of index 2. For any point
  $(x_1, x_2, \cdots, x_d) \in \S$, define $F(x_1, x_2, \cdots, x_d)
  = \prod_{i=1}^{d} f(x_i)$, then it is easy to check that
  $F$ is a Morse function on $s^d\S$. The critical points of $F$
  consist of those $(x_1, x_2, \cdots, x_d)$ where $x_i$ is a critical
  point of $f$,  and the Morse index of $(x_1, x_2, \cdots, x_d)$ is the sum
  of the Morse indices of the $x_i$'s. Therefore,  the number of the critical
  points of $F$ with Morse index $i$ is given by
  $$
  \binom{2g}{i} + \binom{2g}{i-2} + \cdots + \binom{2g}{i-2[i/2]},
  $$
  which is exactly the i-th Betti number of $s^d\S$ (see~\cite{macdonald}).
  Hence, $F(x_1, x_2, \cdots, x_d)$ is a perfect Morse function on $s^d\S$.
  Then we can perturb the Chern-Simons Seiberg-Witten functional such
  that there exists a one-to-one correspondence between the perturbed
  Seiberg-Witten monopoles on $\Y$ and the critical points
  of $F$ on $s^d\S$. Both sets of critical points are non-degenerate and have
  the same relative indices modulo $2|r|$. This implies that
  $\dim HFSW^*(Y,\frs_r) \leq \dim H^*(s^d\S)$.
\hfill $\Box$ \hspace{2mm}

To shorten the notation, we shall write from now on
  \begin{equation}
    V_r=HFSW^*(Y,\frs_r),
  \label{eqn:5}
  \end{equation}
for $-(g-1) \leq r \leq g-1$, $r\neq 0$. In this section we will study the
finite dimensional vector spaces $V_r$ for $r\neq 0$. They have
a natural $\ZZ/2|r|\ZZ$-grading. The only tools we shall
use are the bound on the dimension provided by theorem~\ref{thm:key} and
the gluing theorem~\ref{thm:gluing}. First, it is easily seen that the
diffeomorphism $f\x c: \Y \ar\Y$, where $f:\S \ar \S$ is an orientation
reversing diffeomorphism and $c:\SS^1 \ar \SS^1$ is the conjugation, induces an
isomorphism $V_r \iso V_{-r}$. Henceforth we shall suppose $r>0$
in~\eqref{eqn:5}.

Let $A=\S\x D^2$ be the $4$-manifold given as the product of $\S$
times a $2$-dimensional disc, so that $\bd A=\Y$. Let $\D=\point
\x D^2 \subset A$. The $\spinc$ structures on $A$ are parametrised
by $H^2(A;\ZZ)= H^2(\S;\ZZ) \iso \ZZ$. We write $\frs_r$ for the
$\spinc$ structure on $A$ with $c_1(\frs_r)= 2r \PD[\D]$ (we use
the same name $\frs_r$ for $\spinc$ structures on $Y$ and on $A$.
No confusion should arise from this, as they are compatible in the
sense that $\frs_r|_Y=\frs_r$). Note that $\ker (c_1(\frs_r)) =
\im (H^1(A;\ZZ) \to H^1(Y; \ZZ))$ and $HFSW^*(Y,\frs_r) \cong
HFSW^*_{[\ker (c_1(\frs_r))]}(Y,\frs_r)$, the relative
Seiberg-Witten invariants of $A$ give a map
  \begin{equation}
    \begin{array}{ccc} \AA(\S) & \ar & V_r=HFSW^*(Y,\frs_r) \\
    z & \mapsto & \p_A(\frs_r,z). \end{array}
  \label{eqn:6}
  \end{equation}
As $S= A \cup_Y A =\S \x \SS^2$, the gluing theorem~\ref{thm:gluing} yields
  \begin{equation}
    \sum_{n\in \ZZ} SW_{S,\frs_r +n[\S]}(z_1z_2)
    = \la \p_A(\frs_r,z_1),\p_A(\frs_r,z_2)\ra,
  \label{eqn:7}
  \end{equation}
for any $z_1,z_2 \in \AA(\S)$, where the $\spinc$ structure $\frs_r$ on $S$
is the one with $c_1(\frs_r)= 2r \PD [\SS^2]$. The metric that we must use for
the Seiberg-Witten invariants in the left hand side of~\eqref{eqn:7}
is one giving a long neck, i.e.\ with period point $\om_g$ close to $[\S]$
in $\cK_0=\{ a[\SS^2]+b[\S] /a,b >0\}$. This implies that
$c_1(\frs_r+n[\S]) \cdot \om_g > 0$ as $r>0$, so the invariants are
calculated in the component $-\cK_0$ of the positive cone.

As $r \neq 0$, there is at most one $n\in\ZZ$ that contributes to the
left hand side in~\eqref{eqn:7},
since $c_1(\frs_r +n[\S])=2r [\SS^2]+2n [\S]$ and
$$
  d(\frs_r +n[\S]) = 2rn +2(g-1).
$$
It is thus important to know the Seiberg-Witten invariants of $S=\S \x \SS^2$
for the component $-\cK_0$, which we describe now. We fix the homology orientation
given by the usual orientation of $H^1(S)=H^1(\S)$ and the orientation of
$H^2_+(S)=\RR\om_g$ determined by $-\om_g$.

Fix a symplectic basis $\{\seq{\g}{1}{2g}\}$ of
$H_1(\S;\ZZ)$ with $\g_i \g_{i+g}=\point$, for $1 \leq i \leq g$. Then
 $$
 \AA(\S)=\QQ[ x ] \ox \L^* (\seq{\g}{1}{2g})
 $$
and there is an action of the mapping class group of $\S$,
$\pi_0(\Diff (\S))$, factoring through an action of the symplectic
group $\Spz$, on both $\AA(\S)$ and $V_r$, making the
map~\eqref{eqn:6} equivariant. Define
 $$
 \h= \sum_{i=1}^g \g_i\g_{g+i} \in \L^2 H_1(\S).
 $$
Then the invariant part $\AA(\S)_I$ of $\AA(\S)$ is generated by
$x$ and $\h$. We decompose $\AA(\S)$ in irreducible
$\Spz$-representations as
 $$
  \AA(\S)= \bigoplus_{k=0}^g  \L_0^k \ox {\QQ[x,\h]\over (\h^{g+1-k})},
 $$
where
 $$
 \L_0^k=\L_0^k H_1(\S)=\ker (\h^{g-k+1}: \L^k H_1(\S) \ar \L^{2g-k+2}
 H_1(\S))
 $$
is the primitive component of $\L^k H_1(\S)$, for $0 \leq k \leq
g$. Then, as the Seiberg-Witten invariant $SW_{S,\frs}(z)$ is
invariant under the action of $\Diff(\S)$, $SW_{S,\frs}(z)=0$ for
any $z \in \bigoplus\limits_{k=1}^g  \L_0^k \ox {\QQ[x,\h]/
(\h^{g+1-k})}$, and it only matters to compute $SW_{S,\frs}(z)$
for $z=x^a\h^b$.

\begin{lem}
\label{lem:invariants}
  Fix $0 <r \leq g-1$ and $n\in \ZZ$. Set $d=g-1-r$. Then
  $SW_{S,\frs_r +n[\S]}$ is zero unless $n \leq -1$ and $D=rn+g-1 \geq 0$
  (there is only a finite number of such $n$). In that case
  $SW_{S,\frs_r +n[\S]}(x^a\h^b)={g! \over (g-b)!} (-n)^{g-b}$, for $a+b=D$,
  $0\leq b \leq g$.
  Note that $D \leq d$ and $D\equiv d\pmod r$. As a consequence, for $n=-1$
  (i.e.\ $D=d$) we have $SW_{S,\frs_r - [\S]}(z)=\la z,[s^d \S]\ra$, for any
  $z\in \AA(\S)$ of degree $2d$.
\end{lem}

\noindent {\em Proof.}
  Let $L$ be the determinant bundle of $\frs_r+n[\S]$, so that $c_1(L)=
  2r[\SS^2] +2n[\S]$. Let $H=\S+\f \SS^2$ be a polarisation close to $[\S]$,
  i.e.\ $\f>0$ small. Then $\deg_H L =2r+2n\f >0$, so by~\cite[proposition
  27]{brussee} the non-perturbed Seiberg-Witten moduli space on $S$ is
  $\PP(H^0(K\ox \cL^{\vee})^*)$, where $-K+2\cL=L$, so $K-\cL=\frac{K-L}{2}
  \equiv (g-1-r)[\SS^2] +(-1-n)[\S]$. For $n\geq 0$ this is empty and hence
  $SW_{S,\frs_r +n[\S]}=0$.

  For $n\leq -1$, $d(\frs_r+n[\S])=2 (rn +g-1)$. Let $H_0=\f \S+\SS^2$ be a
  polarisation close to $[\SS^2]$, i.e.\ $\f>0$ small. Then $\deg_{H_0}
  L =2r\f+2n <0$, so by~\cite[proposition
  27]{brussee} the non-perturbed
  Seiberg-Witten moduli space on $S$ is $\PP(H^0(\cL)^*)$, where $\cL=
  \frac{K+L}{2}\equiv (g-1+r)[\SS^2] +(-1+n)[\S]$. Hence the moduli space is
  empty and the Seiberg-Witten invariant for this polarisation, is zero.
  The Seiberg-Witten invariant $SW_{S,\frs_r +n[\S]}$ is obtained via
  wall-crossing from~\cite{OT}. With the notations therein, $u_c \in
  \L^2 H_1(S;\ZZ)$ is given by
  $u_c(\g_i \wedge \g_j) ={1\over 2}\la \g_i \cup \g_j , c_1(L)
  \ra$, i.e. $u_c= n\h$, and
$$
  SW_{S,\frs_r +n[\S]}(x^a\h^b)= \la\h^b {(-u_c)^{g-b}\over (g-b)!},
  [\Jac\, S] \ra = {g! \over (g-b)!} (-n)^{g-b}.
$$
  The sign is as stated as there is a minus sign coming in as we compute
  the invariants in the component $-\cK_0$ and another minus sign because
  we orient $H^2_+$ with $-\om_g$.

  The last statement follows from~\cite{macdonald}.
\hfill $\Box$ \hspace{2mm}

\begin{prop}
\label{prop:isom-no}
  Fix $0<r \leq g-1$ and put $d=g-1-r$.
  Let $z_i\in \AA(\S)$, $i\in I$, be homogeneous elements such that
  $\{ z_i \}_{i\in I}$
  is a basis for $H^*(s^d\S)$, under the epimorphism~\eqref{eqn:00}.
  Consider for each $i\in I$ the element $e_i=\p_A(\frs_r,z_i)\in V_r
  =HFSW^*(Y,\frs_r)$. Then
  $\{ e_i \}_{i\in I}$ is a basis for $V_r$. Therefore $H^*(s^d\S) \ar V_r$,
  $z_i \mapsto e_i$, is a ($\Spz$-equivariant, non-canonical) isomorphism of vector spaces.
\end{prop}

\noindent {\em Proof.}
  Without loss of generality, we may suppose that $\{ z_i \}_{i\in I}$
  is a basis formed by homogeneous elements with non-decreasing degrees.
  The intersection matrix $(\la z_i, z_j\ra)$ is then of the form
  $$
    \left( \begin{array}{ccccc}
    0 &  \cdots & 0& A_0 \\
    0  &  \cdots & A_1 & 0 \\
    \vdots &  \ddots & \vdots & \vdots \\
    A_{2d} &  \cdots &0 &0
    \end{array} \right) \, ,
  $$
  where $A_i$ are the sub-matrices corresponding to the pairing
  $H^i(s^d\S) \otimes H^{2d-i}(s^d\S) \ar \QQ$. So $\det A_i \neq 0$,
  for $0 \leq i \leq 2d$. By the formula~\eqref{eqn:7} and
  lemma~\ref{lem:invariants}, $\la e_i,e_j\ra=0$ if $\deg z_i+\deg z_j > 2d$
  and $\la e_i,e_j\ra=\la z_i,z_j\ra$ if $\deg z_i+\deg z_j = 2d$.
  Therefore the intersection matrix $(\la e_i, e_j\ra)$ is of the form
  $$
    \left( \begin{array}{ccccc}
    \hbox{*} & \cdots & \hbox{*}& A_0 \\
    \hbox{*} & \cdots & A_1 & 0 \\
    \vdots   & \ddots & \vdots & \vdots \\
    A_{2d} & \cdots & 0 & 0
    \end{array} \right) \, ,
  $$
  which is invertible. This implies in particular that
  $\dim V_r \geq \dim H^*(s^d\S)$. As we already have the opposite
  inequality from theorem~\ref{thm:key}, it must be $\dim V_r = \dim H^*(s^d\S)$
  and $\{ e_i \}_{i\in I}$ is a basis for $V_r$.
\hfill $\Box$ \hspace{2mm}

The proof of this proposition shows that the map~\eqref{eqn:6} is
surjective. We have the following

\begin{criterium}
\label{crit:vanishing}
  Let $z\in \AA(\S)$ and $0<|r|\leq g-1$. Then the following are equivalent:
\begin{itemize}
 \item $\p_A(\frs_r,z)=0$.
 \item $SW_{S,\frs_r+n[\S]}(z z_i)=0$ for all $i\in I$ and integer $n$.
 \item $SW_{S,\frs_r+n[\S]}(z z')=0$ for all $z'\in \AA(\S)$ and integer $n$.
\end{itemize}
\end{criterium}

\section{Ring structure of $HFSW^*(\Y,\frs_r)$}
\label{sec:5}

Recall our basic set up. We have the three manifold $Y=\Y$
together with the $\spinc$ structure $\frs_r$ with $c_1(\frs_r)=
2r \PD [\SS^1] \in H^2(Y;\ZZ)$, $0< |r| \leq g-1$, and put $d=
g-1-|r|$. We can define a product on $V_r=HFSW^*(Y,\frs_r)$ as
follows. By criterium~\ref{crit:vanishing},
 $$
 \cI_g=\{z \in \AA(\S) \; | \; \p_A(\frs_r,z)=0 \}
 $$
is an ideal of $\AA(\S)$. So we define an associative and graded
commutative ring structure on $V_r$ by
 $$
 \p_A(\frs_r,z_1)\cdot \p_A(\frs_r,z_2)=\p_A( \frs_r,z_1z_2),
 $$
for any $z_1,z_2\in \AA(\S)$. Therefore $V_r=\AA(\S)/ \cI_g$. This
makes the map~\eqref{eqn:6} an epimorphism of rings.

\begin{lem}\label{lem:mu}
 Let $X_1$ be a $4$-manifold with boundary $\bd X_1=Y$ and let
 $\frs$ be a $\spinc$ structure such that $\frs|_Y=\frs_r$. Then
 for any $z_1\in \AA(X_1)$ and $z_2\in \AA(\S)$ we have
 $$
   \p_A(\frs_r,z_2) \cdot \p_{X_1}(\frs,z_1) = \p_{X_1}(\frs,z_2z_1).
 $$
\end{lem}

\noindent {\em Proof.}
 First, for any $\phi \in V_r$ we have $\la \phi \cdot
 \p_A(\frs_r,z_2),\p_A(\frs_r,z)\ra
 = \la \phi, \p_A(\frs_r,z_2 z) \ra$.
 This is true since by the very definition of the product
 it holds for the elements
 $\phi=\p_A(\frs_r,z')$, which generate $V_r$.

 Now for any $z \in \AA(\S)$ we have
 \begin{eqnarray*}
 & & \la \p_{X_1}(\frs,z_1)\cdot\p_A(\frs_r,z_2), \p_A(\frs_r,z)\ra
 =  \la \p_{X_1}(\frs,z_1),\p_A(\frs_r,z_2 z)\ra= \\
 & & = SW_{X,\frs_r+n[\S]}(z_1z_2z)=
 \la \p_{X_1}(\frs,z_1z_2), \p_A(\frs_r,z)\ra,
 \end{eqnarray*}
 where $X=X_1\cup_Y A$ and $n$ is a suitable integer.
 By criterium \ref{crit:vanishing} we have the result.
\hfill $\Box$ \hspace{2mm}

Note that the isomorphism $V_r \iso V_{-r}$ intertwines the ring
structures, so we may restrict to the case $r>0$. Recall that
$d=g-1-r$.

\begin{thm}
\label{thm:deform}
  Denote by $\cdot$ the product induced in $H^*(s^d\S)$ by the product
  in $V_r$ under the isomorphism of proposition~\ref{prop:isom-no}. Then
  $\cdot$ is a deformation of the cup product graded modulo
  $2r=2(g-1-d)$, i.e.\ for $f_1 \in H^i (s^d\S)$, $f_2 \in H^j(s^d\S)$,
  it is $f_1 \cdot f_2 = \sum_{m \geq 0} \P_m(f_1,f_2)$, where
  $\P_m \in H^{i+j + 2mr}(s^d\S)$ and $\P_0= f_1 \cup f_2$.
\end{thm}

\noindent {\em Proof.}
  By lemma~\ref{lem:invariants}, for any $i, j \in I$,
  $\la e_i,e_j \ra$ is zero unless $\deg z_i + \deg z_j = 2d - 2mr$,
  with $m \geq 0$. Moreover, when $\deg z_i + \deg z_j = 2d$, it is
  $\la e_i, e_j\ra = \la z_i, z_j\ra $.
  Now the same argument as in~\cite[theorem 5]{munoz} accomplishes the
  result, the only difference being that, in our present case, the
  deformation produces terms of increasing degrees.
\hfill $\Box$ \hspace{2mm}

\begin{cor}
\label{cor:vanish}
  Let $f \in \AA(\S)$ be an homogeneous element of degree strictly
  bigger than $2d$. Then $f$ is zero in $V_r$.
\hfill $\Box$
\end{cor}

The last ingredient that we need in order to describe $V_r$ is the
cohomology $H^*(s^d\S)$ of the $d$-th symmetric product $s^d\S$ of
the surface $\S$. Here $d=g-1-r$, so $d$ is in the range $0 \leq d
< g-1$. This cohomology ring was initially described
in~\cite{macdonald} and revisited in~\cite[section 4]{adjuncti}
where it was described as $\Spz$-representation, which is the form
well suited for our purposes. Put an auxiliary complex structure
on $\S$ and interpret $s^d\S$ as the moduli space of degree $d$
effective divisors on $\S$. Let $D \subset s^d\S \x\S$ be the
universal divisor. Then
 $$
 \left\{ \begin{array}{l} \eta= c_1(D) / x \in H^2(s^d\S)
 \\ \q_i=c_1(D) /\g_i \in H^1(s^d\S), \qquad 1\leq i\leq 2g
 \end{array} \right.
 $$
are generators of the ring $H^*(s^d\S)$, i.e.\ there is a graded
$\Spz$-equivariant epimorphism
\begin{equation}
  \AA(\S) \iso \QQ[\eta]\otimes  \L^*(\seq{\q}{1}{2g}) \surj H^*(s^d\S).
\label{eqn:00}
\end{equation}
We set $\h=\sum_{i=1}^g \q_i\q_{g+i} \in H^2(s^d\S)$, abusing a
little bit notation since it correspond to the element $\h$ under
\eqref{eqn:00}. Also we identify $\L^k_0=\L^k_0(\seq{\q}{1}{2g})$
under the same map. Clearly $\eta$ and $\h$ generate the invariant
part $H^*(s^d\S)_I$. The description of $H^*(s^d\S)$ as
$\Spz$-representation is given in the following

\begin{prop}{\bf (\cite[proposition 3.5]{adjuncti})}
\label{prop:hsdS}
  For $0\leq d \leq g-1$ there is a presentation
  $$
    H^*(s^d\S) =\bigoplus_{k=0}^{d} \L^k_0 \ox {\Ceh \over J_k^g },
  $$
  where $J_k^g =(R_k^g ,\h R_{k+1}^g,\h^2 R_{k+2}^g,\ldots, \h^{d+1-k})$,
  $0 \leq k \leq d$, and
  $$
   R_k^g= \sum\limits_{i=0}^{\a} {{(d-k)-\a+1 \choose i} \over{g-k \choose i}}
   {(-\h)^i \over i!}\eta^{\a-i},
  $$
  for $0 \leq k \leq d$, with $\a=[{d-k \over 2}]+1$ (for consistency,
  $R_{d+1}^g=1$). Actually $J_k^g=(R_k^g,\h R_{k+1}^g)$. A basis for $\Ceh /J_k^g$
  as vector space is given by $\eta^a\h^b$, with $2a+b \leq d-k$.
\hfill $\Box$
\end{prop}

For the space $V_r$ we set
 $$
 \left\{ \begin{array}{l} \eta=\p_A(\frs_r,x) \in V_r
 \\ \q_i =\p_A(\frs_r,\g_i) \in V_r, \qquad 1\leq i\leq 2g
 \end{array} \right.
 $$
where $\eta$ has degree $2$ and $\q_i$ degree $1$. (We name with
the same letters elements in $V_r$ and in $H^*(s^d\S)$ as they are
obviously related.) These elements are generators of $V_r$ as
algebra. This means that~\eqref{eqn:6} is a $\Spz$-equivariant
epimorphism
 $$
  \AA(\S) \iso \QQ[\eta]\otimes  \L^*(\seq{\q}{1}{2g}) \surj V_r.
 $$
Clearly $\eta$ and $\h=\sum_{i=1}^g \q_i\q_{g+i}$ generate the
invariant part of $V_r$. Now we are going to relate the ring
structure of $H^*(s^d\S)$ with that of $V_r$.

\begin{prop}
\label{prop:Vr}
  Let $0< r \leq g-1$ and set $d=g-1-r$. Then there is a presentation
  $$
  V_r= \bigoplus_{k=0}^{d} \L^k_0 \ox {\Ceh \over I_k^g},
  $$
  where $I_k^g=(\cR_k^g ,\h \cR_{k+1}^g) \subset \Ceh$ are
  ideals (dependent on $g$, $k$ and $r$) such that
  \begin{equation}
   \cR_k^g= R_k^g + \sum^{\a+mr}_{i=2\a+2mr-(d-k) \atop m>0}
   {a_{im}\over i!{g-k \choose i}} \eta^{\a+mr-i}\h^i,
  \label{eqn:Vr}
  \end{equation}
  where $\a=[{d-k \over 2}]+1$, $R_k^g$ are given in
  proposition~\ref{prop:hsdS}
  and $a_{im}$ are some complex numbers (dependent on $g$, $k$, $r$).
  A basis for $\Ceh / I^g_k$ is given by $\eta^a\h^b$, with $2a+b \leq d-k$.
\end{prop}

\noindent {\em Proof.}
  Let $\{z_i^{(k)}\}$ be a basis for $\L^k_0$. Then by proposition~\ref{prop:hsdS},
  $z_i^{(k)} x^a\h^b$, $2a+b+k \leq d$, form a basis for $H^*(s^d\S)$.
  We use this basis in proposition~\ref{prop:isom-no} to construct a
  ($\Spz$-equivariant) isomorphism $H^*(s^d\S) \iso V_r$.
  The fact that $R_k^g \in J_k^g$ means that $z_i^{(k)} R_k^g=0$ in
  $H^*(s^d\S)$. Fix $z_0^{(k)}=\q_1\cdots\q_k \in \L^k_0$, then
  $ \L^k_0 = \text{Span} <\Spz z_0^{(k)}>$.
  Rewriting $z_0^{(k)} R_k^g$ in terms of the product $\cdot$
  of theorem~\ref{thm:deform}, and using the arguments
  of~\cite[section 2]{Intri} (and the fact that the action of $\Spz$
  is compatible with the ring structure on $V_r$), we get that
  \begin{equation}
    z_0^{(k)} R_k^g + \sum_{m>0} z_i^{(k)} R_{kim}^g =0
  \label{eqn:zi}
  \end{equation}
  in $V_r$, where $\deg R_{kim}^g=\deg R_k^g +mr=\a+mr$,
  and $R_{kim}^g$ is expressible
  in terms of the chosen basis, i.e.\ as a linear combination of the
  monomials $\eta^{\a+mr-j}\h^j$, for $2\a+2mr-(d-k)\leq j\leq\a+mr$.
  As in the proof
  of~\cite[proposition 16]{munoz2}, we have that the only nonvanishing
  $R_{kim}^g$ in~\eqref{eqn:zi} correspond to $z_0^{(k)}$
  (otherwise one can find an element of $\Spz$ only fixing $z_0^{(k)}$,
  which
  would produce a relation between the elements of the basis of $V_r$, which
  is impossible), so~\eqref{eqn:zi}
  reduces to $z_0^{(k)} (R_k^g + \sum_{m>0} R_{k0m}^g) =0$ in $V_r$.
  This produces the relation $\cR_k^g=R_k^g + \sum_{m>0} R_{k0m}^g$
  as stated in~\eqref{eqn:Vr}.

  Also $\h \cR_{k+1}^g\in I_k^g$ since
  $\h I_{k+1}^g \subset I_k^g$. Now $I_k^g$ is generated by these two
  elements since $J_k^g$ is generated by $R^g_k$ and $\h R^g_{k+1}$
  (see~\cite[section 2]{Intri}).
\hfill $\Box$ \hspace{2mm}

\begin{rem}
\label{rem:d-k}
  Note that for $d$ odd, $\cR_k^g$ is the relation uniquely determined by
  expressing $\eta^{\a} \in \Ceh/I^g_k$, $\a=[{d-k \over 2}]+1$,
  in terms of the monomials of the
  basis $\eta^a\h^b$, $2a+b \leq d-k$.
  For $d$ even, $\cR_k^g$ is the relation uniquely determined by
  expressing $\eta^{\a} -{(d-k)-\a+1 \over g-k}\eta^{\a-1}\h$
  in terms of the monomials of the basis.
\end{rem}

\begin{cor}
\label{cor:hom-symm}
  There is an isomorphism of associated graded rings
  $$\Gr_{\h} V_r \iso \Gr_{\h} H^*(s^d \S),$$
  where $d=g-1-|r|$. Let $HF^*(\Y)$ be the instanton Floer homology
  of $Y=\Y$ with $SO(3)$-bundle with $w_2=\PD[\SS^1]$, which was
  computed in~\cite{munoz}. This can be decomposed~\cite[proposition
  20]{munoz} according to the eigenvalues of $\a=2\mu(\S)$ as
  $HF^*(\Y)=\oplus_{r=-(g-1)}^{g-1} H_r$, where $\a$ has eigenvalue
  $4r$ (if $r$ is odd) or $4\ima r$ (if $r$ is even) on $H_r$,
  $-(g-1) \leq r\leq g-1$. Then~\cite[corollary 3.7]{adjuncti}
  gives an isomorphism
  $$\Gr_{\h} V_r \iso \Gr_{\g} H_r,$$
  where $\g=-2\sum \mu(\g_i)\mu(\g_{g+i})$.
\hfill $\Box$
\end{cor}

\begin{lem}
\label{lem:Ikg}
  Fix $r>0$. Then the ideals of proposition~\ref{prop:Vr}
  satisfy the recursion $I^g_k=I^{g-1}_{k-1}$, for $k>0$ and $r\leq g-2$.
  Equivalently, $\cR_k^g=\cR^{g-k}_0$.
\end{lem}

\noindent {\em Proof.}
  By the computation of the Seiberg-Witten invariants
  of $S$ in lemma~\ref{lem:invariants} and the invariance under the action
  of $\Spz$, we have
\begin{eqnarray*}
  SW_{S,\frs_r +n[\S]}(\g_1\cdots \g_k \g_{g+1}\cdots\g_{g+k} x^a\h^b)
  &=& SW_{S,\frs_r +n[\S]}\left({1 \over k!} {g \choose k}^{-1} \h^k x^a\h^b\right)\\
  &=& {(g-k)! \over (g-k-b)!} (-n)^{g-k-b},
\end{eqnarray*}
  for $a+b =g-k -1-rn$. Therefore for any $R\in \AA(\S)_I$, $z\in \AA(\S)$,
  $$
  \la \q_1\cdots \q_{k-1} \q_g R, \q_{2g} z \ra_g=
  \la \q_1\cdots \q_{k-1} R, z \ra_{g-1},
  $$
  where the subindex means the genus of the surface $\S$.
  This implies the statement by criterium \ref{crit:vanishing}.
  The last part follows from remark \ref{rem:d-k}.
\hfill $\Box$ \hspace{2mm}

Now we aim to compute the coefficients $a_{im}$ of $\cR_0=\cR_0^g$
in~\eqref{eqn:Vr}.
Let $m>0$. We collect the coefficients together in a polynomial
\begin{equation}
\label{eqn:0}
    p_m(x)=\sum a_{im} x^{g-i},
\end{equation}
where we consider $a_{im}=0$ for $i \notin [2\a+2mr-d, \a+mr]\cap
\ZZ$, $\a=[{d\over 2}]+1$. Note that there are a finite number of
non-zero polynomials. By analogy we consider
\begin{equation}
\label{eqn:1}
  p_0(x)=(x-1)^{d-\a+1}x^{g-(d-\a-1)},
\end{equation}
so that $R_0 =\sum {a_{i0}\over i!{g \choose i}} \eta^{\a-i}\h^i$, as given
in proposition~\ref{prop:hsdS}. By definition $\cR_0=0 \in V_r$, therefore we
have $\la \cR_0, \eta^a\h^b\ra=0$, whenever $\a+a+b=d-kr$, $k\geq 0$.
Now using the computation of the invariants of $S$ in
lemma~\ref{lem:invariants},
we get
$$
  \sum_{m=0}^k {a_{im}\over i!{g \choose i}} {g! \over (g-b-i)!} (k-m+1)^{g-b-i}=
  \sum_{m=0}^k a_{im} {(g-i)! \over (g-b-i)!}(k-m+1)^{g-b-i}=0,
$$
for all $k\geq 0$ and $0\leq b \leq d-\a-kr$. So
$$
  \sum_{m=0}^k {d^b \over dx^b} p_m(x) \Big|_{x=k-m+1} =0.
$$
for all $k\geq 0$ and $0\leq b \leq d-\a-kr$.
 By Taylor expansion, this is equivalent to saying that
\begin{equation}
\label{eqn:2}
  p_k(x) \equiv -\left( p_0(x+k) +p_1(x+k-1) + \cdots +p_{k-1}(x+1) \right)
  \pmod{(x-1)^{d-\a-kr+1}}.
\end{equation}
This condition, together with the fact that $p_k(x)$ has degree
$g-(2\a+2kr-d)$ and it
is divisible by $x^{g-(\a+kr)}$, uniquely determines $p_k(x)$ by recursion.

For instance, let us calculate explicitly $p_1(x)$. From~\eqref{eqn:1}
we have that
\begin{eqnarray*}
  p_0(x+1)&=&x^{d-\a+1}(x+1)^{g-(d-\a-1)}= x^{g-\a-r}(x+1)^{\a+r}\\
  &=& x^{g-\a-r} \sum_{k=0}^{\a+r} {\a+r \choose k} 2^k (x-1)^{\a+r-k},
\end{eqnarray*}
using that $d=g-1-r$. Now $p_1(x)$ is divisible by $x^{g-\a-r}$, has degree
$(g-\a-r)+(d-\a-r)$ and $p_1(x) \equiv -p_0(x+1) \pmod{(x-1)^{d-\a-r+1}}$.
Therefore
\begin{eqnarray*}
  p_1(x) & = & - x^{g-\a-r}  \sum_{k=2\a+2r-d}^{\a+r}
  \hspace{-2mm} {\a+r \choose k} 2^k (x-1)^{\a+r-k}= \\
   &= &- x^{g-\a-r} \hspace{-2mm} \sum_{2\a+2r-d \leq k\leq \a+r\atop 0\leq j\leq \a+r-k}
  \hspace{-2mm} 2^k {\a+r \choose k} {\a+r-k\choose j}(-1)^j x^{\a+r-k-j}= \\
  &=& \hspace{-2mm} \sum_{2\a+2r-d \leq k\leq \a+r\atop 0\leq j\leq \a+r-k}
  \hspace{-3mm} (-1)^{j+1} {(\a+r)! \over k!j!(\a+r-k-j)!} 2^k x^{g-k-j}.
\end{eqnarray*}
 From this we may write the coefficients $a_{i1}$ as
$$
  a_{i1}= \sum_{j=0}^{i-(2\a+2r-d)} (-1)^{j+1}
  {(\a+r)! \over (i-j)!j!(\a+r-i)!} 2^{i-j},
$$
for $2\a+2r-d \leq i\leq \a+r$.

We can compute the rest of the coefficients $a_{im}$, for $m>1$,
by recurrence using this method but the result is a collection of rather
cumbersome formulae which do not shed light on the ring structure of $V_r$.
This is to no surprise: the shape of the relations $\cR_k^g$ depends on the
basis of $\Ceh/I_k^g$ that we have chosen in proposition~\ref{prop:Vr}, and
this basis has been chosen rather arbitrarily. We shall present now a
slightly modified version of the previous argument which computes explicitly
(a full set of) relations for $V_r$, by just not fixing any basis for
$\Ceh/I_k^g$. This leads to a closed formula for generators of the ideals $I_k^g$.

\begin{thm}
\label{thm:cRk}
  Let $0< r \leq g-1$ and set $d=g-1-r$. Then there is a presentation
  $$
  V_r= \bigoplus_{k=0}^{d} \L^k_0 \ox {\Ceh \over
  (\tilde\cR_k^g ,\h \tilde\cR_{k+1}^g, \eta^{d+1},\h^{d+1})},
  $$
  where
  $$
   \tilde\cR_k^g= \sum_{i=0}^{\a}
   {{d-k-\a+1 \choose i}\over i!{g-k \choose i}}(-1)^i \eta^{\a-i}\h^i -
   \sum_{i=0}^{\a+r} {{\a+r\choose i}\over i!{g-k \choose i}} \eta^{\a+r-i}\h^i,
  $$
  where $\a=[{d-k \over 2}]+1$, for $0 \leq k \leq d$, and $\tilde\cR_{d+1}^g=1$.
\end{thm}

\noindent {\em Proof.}
  By lemma~\ref{lem:Ikg} it is enough to find a relation for $k=0$,
  \begin{equation}
   \tilde\cR_0= R_0 + \sum_{m>0 \atop 0 \leq i\leq\a+mr}
   {a_{im}\over i!{g \choose i}} \eta^{\a+mr-i}\h^i,
  \label{eqn:cRk}
  \end{equation}
  This time we do not restrict the range for $i$. We only note
  that we can suppose $a_{im}=0$ if $i> g$, since $\h^{g+1}=0$.
  As before, we collect the coefficients $a_{im}$ of~\eqref{eqn:cRk}
  in a polynomial $p_m(x)=\sum a_{im} x^{g-i}$, where $a_{im}=0$
  for $i \notin [0, \a+mr]\cap \ZZ$. Also $p_0(x)=(x-1)^{d-\a+1}x^{g-(d-\a-1)}$.
  The condition that $\cR_0$ be a relation is translated into
\begin{equation}
\label{eqn:20}
  p_k(x) \equiv -\left( p_0(x+k) +p_1(x+k-1) + \cdots +p_{k-1}(x+1) \right)
  \pmod{(x-1)^{d-\a-kr+1}}.
\end{equation}
  We want to find polynomials $p_k(x)$ of degree $g$ solving~\eqref{eqn:20}.
  This time the $p_k(x)$ are not determined uniquely, but we only need
  to find one solution. Since $p_0(x+1)=x^{d-\a+1}(x+1)^{g-(d-\a-1)}=
  x^{g-\a-r}(x+1)^{\a+r}$, we may choose
$$
  p_1(x)= -x^{g-\a-r}(x+1)^{\a+r}
$$
  and $p_k(x)=0$ for $k\geq 2$. This gives $a_{i1}= -{\a+r \choose i}$,
  $0 \leq i\leq \a+r$, and $a_{im}=0$ for $m\geq 2$. In this way we have
  found $\tilde\cR^g_k, \h \tilde \cR^g_{k+1} \in I_k^g$ as given in the
  statement. However they do not
  generate the whole ideal as may be seen by looking at the associated
  graded ring $\Gr_{\h} \left( \Ceh/ (\tilde\cR^g_k, \h \tilde \cR^g_{k+1})
  \right)$, so we need to add more relations.
  The nilpotence relations $\eta^{d+1},\h^{d+1}$ are always
  satisfied by corollary~\ref{cor:vanish}. To see that these relations
  suffice, write any $f \in I_k^g$ as $f=a_1 \cR^g_k+a_2\h \cR^g_{k+1}$,
  by proposition~\ref{prop:Vr}.
  Then $f- a_1\tilde\cR^g_k-a_2\h \tilde \cR^g_{k+1} \in I_k^g$ and has
  higher degree than that of $f$. Proceed recursively until we get
  a polynomial in $( \eta^{d+1},\h^{d+1})$.
\hfill $\Box$ \hspace{2mm}

\section{Seiberg-Witten invariants of connected sums along surfaces}
\label{sec:7}

We want to show, as a first application, how the knowledge of the
previous sections can be used to compute the Seiberg-Witten
invariants of $4$-manifolds which appear as connected sums along
surfaces of other $4$-manifolds. This was first dealt with in a
particular case in~\cite{taubes} to get a proof of the symplectic
Thom conjecture. In the context of Donaldson invariants it has
been extensively treated in~\cite{genusg}~\cite{hff}.

The set up is as follows (see~\cite{genusg}). Let $\bar X_1$ and
$\bar X_2$ be smooth oriented $4$-manifolds and let $\S$ be a
compact oriented surface of genus $g\geq 2$. Suppose that we have
embeddings $\S \inc \bar X_i$ with image $\S_i$ representing a
{\em non-torsion} element in homology whose self-intersection is
zero. This implies that $b^+>0$. Now take small closed tubular
neighbourhoods $N_{\S_i}$ of $\S_i$ which are isomorphic to $A=
\S\x D^2$. Let $X_i$ be the closure of $\bar X_i - N_{\S_i}$,
$i=1,2$. Then $X_i$ is a $4$-manifold with boundary $\bd X_i=Y=\Y$
and $\bar X_i=X_i \cup_Y A$. Take an identification $\phi: \bd X_1
\ar \, -\bd X_2$ (i.e. an orientation reversing bundle
isomorphism). We define the connected sum of $\bar X_1$ and $\bar
X_2$ along $\S$ as
 $$
  X=X(\phi) = X_1 \cup_{\phi} X_2.
 $$
The resulting $4$-manifold depends in general on the isotopy class
of $\phi$, but we shall drop $\phi$ from the notation when there
is no danger of confusion, and write then $X=\bar X_1 \#_{\S} \bar
X_2$. Consider $\spinc$ structures $\frs_i$ on $X_i$ such that
$\frs_1|_Y \iso -\frs_2|_Y \iso \frs_Y$, with $c_1(\frs_Y)\neq 0$,
so that they can be glued together to get a $\spinc$ structure
$\frs_o$ on $X$. The $\spinc$ structures $\frs$ such that
$\frs|_{X_i}=\frs_i$, $i=1,2$, are those of the form $\frs_o+h$,
where $h$ is an element in the image of  $H_2(Y;\ZZ) \ar
H_2(X;\ZZ)\iso H^2(X;\ZZ)$, where the last map is Poincar\'e
duality. Let $\rim \subset H^2(X;\ZZ)$ be the subspace generated
by the rim tori~\cite{rim}, i.e.\ the image of $H_1(\S;\ZZ)\ox
[\SS^1] \subset H_2(Y;\ZZ)$ in $H^2(X;\ZZ)$. Then any $\spinc$
structure $\frs$ such that $\frs|_{X_i}=\frs_i$, $i=1,2$, is of
the form $\frs=\frs_o+h+n\S$, where $h\in \rim$, $n\in \ZZ$.

If $\frs_Y$ has $c_1(\frs_Y) \neq 2r\PD [\SS^1]$, for any
$-(g-1)\leq r \leq g-1$, $r\neq 0$, then theorem~\ref{thm:key}
tells us that $SW_{X,\frs}=0$. Now consider the case
$\frs_Y=\frs_r$, with $-(g-1)\leq r \leq g-1$, $r\neq 0$. Set
$d=g-1-|r|$ as usual.

\begin{thm}
\label{thm:gl-neq}
  Fix $z_i \in \AA(\S)$, $i \in I$, homogeneous elements
  such that $\{ z_i \}_{i \in I}$ is a
  basis for $H^*(s^{d}\S)$. Then there exists a universal matrix
  $(m_{ij})_{i,j \in I}$ such that for every connected sum
  $X=\bar X_1 \#_{\S} \bar X_2$ along a surface of genus $g$, $\spinc$
  structures $\bar{\frs}_i$ on $\bar X_i$ with $c_1(\bar{\frs}_i)
  \cdot \S =2r$ and $\spinc$ structure $\frs_o$ on $X$ obtained by
  gluing $\bar\frs_1$ and $\bar\frs_2$, we have
  $$
    \sum_{h\in\rim} SW_{X,\frs_o+h}
   (z_1z_2)= \sum_{n,m \in \ZZ \atop i,j\in I}
   m_{ij}  SW_{\bar X_1,\bar{\frs}_1+ n\S}(z_1 z_i)\cdot
   SW_{\bar X_2,\bar{\frs}_2+m\S}(z_2z_j) ,
  $$
  for any $z_1 \in \AA(\bar X_1)$ and $z_2 \in \AA(\bar X_2)$ with
  $d(\frs_o)=\deg z_1+\deg z_2$ (note that at most one $n$ and one $m$
  appear in every summand of the right hand side).
  If any of the manifolds involved has $b^+=1$ then its
  Seiberg-Witten invariants are computed
  for the component of the positive cone containing $-r\PD[\S]$.
\end{thm}

\noindent {\em Proof.}
  Let $\frs_i=\bar\frs_i|_{X_i}$, $i=1,2$.
  By proposition~\ref{prop:isom-no}, the elements $e_i=\p_A(\frs_r,z_i)$,
  $i\in I$, form a basis for $V_r=HFSW^*(Y,\frs_r)$.
  Therefore $V_r \ar \RR^{|I|}$, given as
  $\phi \mapsto (\la \phi,\p_A(\frs_r,z_i) \ra)_{i\in I}$, is an
  isomorphism such that
  $$
   \p_{X_1}(\frs_1, z_1) \mapsto \left( \sum_{n\in \ZZ} SW_{\bar
   X_1,\bar\frs_1+ n\S} (z_1 z_i) \right)_{i\in I}.
  $$
  Theorem~\ref{thm:gluing} says that
  $$
  \sum_n SW_{\bar X_1, \bar\frs_1+ n\S}(z_1 z_i) =
  \la \p_{X_1}(\frs_1,z_1), e_i\ra, \quad
  \sum_m SW_{\bar X_2, \bar\frs_2+ m\S}(z_2 z_j) =
  \la \p_{X_2}(\frs_2,z_2), e_j\ra
  $$
  and
  $$
    \sum_{\{\frs / \frs|_{X_i}=\frs_i,\; i=1,2 \}}\hspace{-5mm} SW_{X,\frs}
   (z_1z_2)= \la \p_{X_1}(\frs_1,z_1), \p_{X_2}(\frs_2,z_2)\ra.
  $$
  Only the $\spinc$ structures of the form $\frs=\frs_o+h$, $h\in \rim$,
  satisfy $d(\frs)=\deg z_1+\deg z_2$.
  The result follows with
  $(m_{ij})$ being the inverse of the intersection matrix for the basis
  $\{e_i\}_{i\in I}$. Note that this matrix is explicitly computable,
  since by lemma~\ref{lem:invariants} the products $\la e_i,e_j\ra$
  are known.
\hfill $\Box$ \hspace{2mm}

\begin{cor}
\label{cor:st}
  If either of $\bar X_i$ has simple type then
  $\sum_{h \in\rim} SW_{X,\frs_o+h} (x z)=0$, for any $z \in \AA(X)$ and
  $\spinc$ structure $\frs_o$ on $X$ with $c_1(\frs_o)\cdot \S\neq 0$.
  Analogously, if either of $\bar X_i$ has strong simple type then
  $\sum_{h \in\rim} SW_{X,\frs_o+h} (z)=0$, for any $z \in \AA(X)$ with
  $\deg(z)>0$ and any $\spinc$ structure $\frs_o$ on $X$
  with $c_1(\frs_o)\cdot \S\neq 0$.
\hfill $\Box$
\end{cor}

In order to remove the summation over the subspace $\rim$ in
corollary~\ref{cor:st} we need an extra condition.

\begin{defn}
\label{def:admissible}
  A connected sum $X=\bar X_1\#_{\S} \bar X_2$ is admissible if $\rim$ has
  no torsion and there exists
  a subspace $V\subset H^2(X;\ZZ)$ such that $H^2(X;\ZZ)=V\oplus \rim$ and
  $c_1(\frs)\in V$ for every basic class $\frs$ of $X$.
\end{defn}

\begin{rem}
\label{rem:algebraic}
  Suppose that $\bar X_1$ and $\bar X_2$ are K\"ahler surfaces and
  $\S_i \subset \bar X_i$ are smooth complex curves of genus $g$, isomorphic
  as complex curves, such that there is a deformation K\"ahler family $\cZ
  \stackrel{\pi}{\ar} D^2 \subset \CC$ with fiber $Z_t=\pi^{-1}(t)$,
  $t \neq 0$, smooth and $Z_0=\pi^{-1}(0)=\bar X_1\cup_{\S} \bar X_2$, the
  union of $\bar X_1$ and $\bar X_2$ along $\S_1 =\S_2$ with a normal
  crossing.
  Then the general fiber $X=Z_t$ is the connected sum $X=\bar X_1 \#_{\S}
  \bar X_2$ with identification given by the isomorphism between
  the normal bundles
  of $\S_1$ and $\S_2$. If $H^2(X;\ZZ)$ has no torsion
  then this identification is admissible, since
  for any basic class $\frs$ one has $c_1(\frs) \in H^{1,1}$ and this
  space is orthogonal to $\rim$, as for any $T \in \rim$, it is
  $T^2=0$ and $\om \cdot T=0$ ($\om$ standing for the K\"ahler form). This
  implies that $T \notin H^{1,1}$ unless $T=0$.
\end{rem}

\begin{rem}
\label{rem:vanishing-discs}
  In~\cite[definition 4.1]{MSz}, Morgan and Szab\'o define admissible
  identification when there exists a collection of primitive embedded
  $(-2)$-spheres in $X$ (obtained by pasting embedded $(-1)$-discs in
  $X_1$ and $X_2$) generating a subspace $V \subset H_2(X)$ such
  that $H_2(X)=\cH \oplus V$, where $\cH=\{ D\in H_2(X) / D|_Y =k[\SS^1],
  \text{some $k$}\}$. Then $c_1(\frs)\cdot V=0$ for any
  basic class $\frs$, and this
  implies admissibility in the sense of definition~\ref{def:admissible}
  (assuming again that $H^2(X;\ZZ)$ has no torsion).
\end{rem}

\begin{cor}
\label{cor:st-adm}
  Suppose that the connected sum $X=\bar X_1\#\bar X_2$ is admissible.
  If either of $\bar X_i$ has simple type then $SW_{X,\frs}(xz)=0$ for
  any $z\in \AA(X)$ and any $\spinc$ structure $\frs$ with
  $c_1(\frs)\cdot \S\neq 0$.

  If either of $\bar X_i$ has strong simple type then $SW_{X,\frs}(z)=0$ for
  any $z\in \AA(X)$ with $\deg(z)>0$ and any $\spinc$ structure $\frs$ with
  $c_1(\frs)\cdot \S\neq 0$.
\hfill $\Box$
\end{cor}

The formula in theorem~\ref{thm:gl-neq} becomes simpler when both
$\bar X_i$ have $b_1=0$. We have the following result

\begin{thm}
\label{thm:b1=0}
  Let $X=\bar X_1 \#_{\S} \bar X_2$ be a connected sum
  along a surface of genus $g$ where both
  $\bar X_1$ and $\bar X_2$ are of $H_1$-simple type. Let
  $\bar{\frs}_i$ be $\spinc$ structures on $\bar X_i$ with $c_1(\bar{\frs}_i)
  \cdot \S =2r\neq 0$ and let $\frs_o$ be a $\spinc$ structure on $X$ obtained by
  gluing $\bar\frs_1$ and $\bar\frs_2$. Suppose $d= g-1-|r| \geq 0$.
  Let $z\in \AA(X)$ with $\deg(z)=d(\frs_o)$. Then
  $\sum_{h\in\rim} SW_{X,\frs_o+h} (z)=$
  $$
  \left\{ \begin{array}{ll}
  \sum\limits_{n,m \in \ZZ} (-1)^{d/2} {g-1 \choose d/2}
  SW_{\bar X_1,\bar\frs_1+ n\S}(x^{d/2}) \cdot
  SW_{\bar X_2,\bar\frs_2+m\S} (x^{d/2}), & \text{$z=1$, $d$ even}\\
  0, & \text{otherwise.}
  \end{array}\right.
 $$
  Note that at most one $n$ and one $m$ contribute to this formula.
\end{thm}

\noindent {\em Proof.} By lemma~\ref{lem:mu}, $\q_j
\p_{X_1}(\frs_1,z_1) = \p_{X_1}(\frs_1,\g_j z_1)=0$, for
$j=1,\ldots, 2g$. Therefore $$ \p_{X_1}(\frs_1,z_1) \in K=
\bigcap_{1 \leq j\leq 2g} \ker \q_j, $$ for any $z_1 \in
\AA(X_1)$. Since $K\subset V_r$ is invariant under the action of
$\Spz$, we have that $K \subset \Ceh/I^g_0$ by
proposition~\ref{prop:Vr}. Now $f\in K$ if and only if $\q_j f=0$,
for $j=1,\ldots, 2g$. This means that $f \in I^g_1$ in the
notation of proposition~\ref{prop:Vr}. So $K=I^g_1/I^g_0$, where
the generators of $I^g_1$ are $\cR^g_1$ and $\h \cR^g_2$. The
intersection pairing $\la,\ra: K\ox K \ar \QQ$ is the restriction
of the pairing of $V_r$. Now by lemma~\ref{lem:invariants} $\la
e_i,e_j\ra=0$ if $\deg z_i + \deg z_j >2d$. For $d$ odd, all the
homogeneous components of all the elements in $I^g_1$ have degree
strictly bigger than $d$ (note that the component $R^g_1$ of
$\cR^g_1$ has degree $2([{d-1 \over 2}]+1) =d+1$ and it is the
component of lowest degree). So $K \ox K \ar \QQ$ is the zero map
for $d$ odd, which proves the second line.

For $d$ even, all the homogeneous components of all the
elements in $I^g_1$ have degree strictly bigger than $d$, except for
$R^g_1$, which has degree $d$. So $K \ox K \ar \QQ$ has rank $1$.
Hence, for $z_1\in\AA(\S)$ and $z_2\in \AA(\S)$, we have that
 $$
  \sum_{h\in\rim} SW_{X,\frs_o+h} (z_1z_2)=
  \sum_{n,m\in \ZZ} c\, SW_{\bar X_1,\bar\frs_1+n\S} (z_1x^{d/2})
   SW_{\bar X_2,\bar\frs_2+m\S} (z_2x^{d/2}),
 $$
where we have used $SW_{\bar X_i,\bar\frs_i+n\S} (z_i R^g_1)=
SW_{\bar X_i,\bar\frs_i+n\S} (z_i x^{d/2})$, and with $c=\la
R_g^1,R_g^1 \ra^{-1}$. To compute $c$, note that $\h R_g^1=0$ so
 $$
  \la R_g^1,R_g^1 \ra= \la R_g^1,\eta^{d/2} \ra
  =\sum_{i=0}^{\a} {\binom{d/2}{i}\over i!\binom{g-1}{i}}
   (-1)^i {g!\over (g-i)!} = \sum_{i=0}^{\a} (-1)^i \binom{\a}{i} {g\over
   g-i}= (-1)^{\a} \binom{g-1}{\a}^{-1},
 $$
with $\a=d/2$. Finally corollary~\ref{cor:vanish} implies that
$SW_{\bar X_i,\bar\frs_i+n\S} (z_ix^{d/2})=0$ for any $z_i$ with
$\deg (z_i) >0$. \hfill $\Box$ \hspace{2mm}

The following corollary is analogue to the result
in~\cite[corollary 15]{genusg} regarding the Kronheimer-Mrowka
basic classes.

\begin{cor}
\label{cor:b1=0}
  Suppose $\bar X_1$ is of strong simple type and $\bar X_2$ has
  $H_1$-simple type. Suppose that the connected sum $X=\bar X_1
  \#_{\S}\bar X_2$ is admissible. Then there are no basic classes $\frs$
  with $0<|c_1(\frs) \cdot \S |<2g-2$. The
  basic classes for $X$ with $c_1(\frs)=\pm (2g-2)$
  are indexed by pairs of basic classes
  $(\bar\frs_1, \bar\frs_2)$ for $\bar X_1$ and $\bar X_2$ respectively,
  such that $\bar\frs_1\cdot \S =\bar\frs_2 \cdot \S = \pm (2g-2)$.
\hfill $\Box$
\end{cor}

\section{Higher type adjunction inequalities}
\label{sec:8}

In this section we shall reprove the higher type adjunction
inequalities for non-simple type $4$-manifolds obtained by
Oszv\'ath and Szab\'o in~\cite{OS}. Our method of proof is
considerable simpler and parallels the proof of the higher type
adjunction inequalities in the context of Donaldson invariants
given in~\cite{adjuncti}. On the weak side, we cannot deal with
the case $\S^2=0$, $c_1(\frs)\cdot \S=0$, due to the fact that we
have to restrict to non-torsion $\spinc$ structures for our study
of the Seiberg-Witten-Floer homology.

\vspace{4mm}
\noindent {\em Proof of theorem~\ref{thm:main4}.\/}
  Without loss of generality, by reversing the orientation of $\S$ in the
  case $b^+>1$, we can suppose that $c_1(\frs)\cdot \S \leq 0$. We reduce
  to the case of self-intersection zero by blowing-up. Let $N=\S^2$ and
  consider the blow-up $\tilde X=X\# N \overline{\CP}^2$ with
  exceptional divisors $E_1, \ldots, E_N$. Let $\tilde \S=\S -E_1-
  \cdots -E_N$ be the proper transform of $\S$, which is an embedded
  surface of self-intersection zero and genus $g$, with $b \in \AA
  (\tilde \S)\iso \AA (\S)$.
  Consider the $\spinc$ structure $\tilde\frs$ on $\tilde X$
  with $c_1(\tilde\frs)=c_1(\frs) -E_1 -\cdots -E_N$. Then
  $d (\frs) = d (\tilde\frs)$,
$$
  - c_1(\tilde\frs)\cdot \tilde\S + \tilde\S^2= - c_1(\frs)\cdot \S + \S^2,
$$
  and $SW_{\tilde X,\tilde\frs}(a\, b)=SW_{X,\frs}(a\, b)\neq 0$.

  Therefore we can suppose that $\S^2=0$ and $c_1(\frs)\cdot \S=-2r$, with
  $0< r \leq g-1$. Let $\{\g_i\}$ be a symplectic basis of $H_1(\S)$
  with $\g_i\cdot\g_{g+i}=1$, $1\leq i \leq g$.
  Without loss of generality we may also suppose that  $b=x^p \g_{i_1}
  \cdots \g_{i_m}$, $\deg (b)=2p+m$.
  Now let $A=\S\x D^2$ be a small tubular neighbourhood
  of $\S \subset X$ and consider the splitting $X=X_1 \cup_Y A$, where $X_1$
  is the closure of the complement of $A$ and $\bd X_1=\bd A=Y=\Y$.
  In this case $\frs$ is the only $\spinc$ structure appearing in the
  gluing formula in theorem~\ref{thm:gluing}, so
$$
 0 \neq SW_{X,\frs}(a\, b) = \la \p_{X_1}(\frs, a),\p_A(\frs_{-r},b) \ra.
$$
  (In the case $b^+=1$, the metric giving a long neck for $X$ has
  period point close to $[\S]$. Therefore we are calculating the
  invariants in the component $\cK(X)$ of the positive cone
  containing $\PD[\S]\in H^2(X;\ZZ)$, since $c_1(\frs)\cdot
  \S<0$.)
  Then $\p_A(\frs_{-r},b) \in HFSW^*(Y,\frs_{-r})=V_{-r} \iso V_r$ is
  non-zero and therefore $\eta^p \q_{i_1}
  \cdots \q_{i_m} \neq 0 \in V_r$. By corollary~\ref{cor:vanish},
  this implies
  $2p+m \leq 2d=2(g-1-|r|)$. Therefore $2r +\deg (b) \leq 2g-2$.
\hfill $\Box$

\vspace{4mm}
\noindent {\em Proof of theorem~\ref{thm:main6}.\/}
  Again we may suppose that $\S^2$=0 and $c_1(\frs)\cdot \S=-2r$, with
  $0<r \leq g-1$. Suppose also that  $b=x^p \g_{i_1}
  \cdots \g_{i_m}$, $\deg (b)=2p+m$.
  Now let $A=\S\x D^2$ be a small tubular neighbourhood
  of $\S \subset X$ and consider the splitting $X=X_1 \cup_Y A$. Then
$$
  0 \neq SW_{X,\frs}(a\, b) = \la \p_{X_1}(\frs, a), \p_A(\frs_{-r},b) \ra.
$$
  Here $\p_{X_1}(\frs, a)\in V_{-r}$ lives in the kernels of $\seq{\q}{1}{l}$,
  since as $\imath_*(\g_j)=0 \in H_1(X)$,
  $$
  \q_j\p_{X_1}(\frs, a)= \p_{X_1}(\frs, \g_j a)=0, \qquad j=1,\ldots, l.
  $$
  Therefore it must be $\p_A(\frs_{-r},b)=
  \eta^p \q_{i_1} \cdots \q_{i_m} \notin (\seq{\q}{1}{l})$ in $V_{-r}$.
  The argument in~\cite[proposition 4.5]{adjuncti}
  (using proposition~\ref{prop:Vr})
  shows that any element of degree bigger strictly bigger
  than $g-1-|r|$ must lie in the ideal
  $(\seq{\q}{1}{g-1-|r|})$ of $V_{-r}$. So if $l \geq g-1-|r|$
  then $2p+m \leq g-1-|r|$, i.e.\
  $\deg(b)\leq g-1-|r|$ and $|2r| +2\deg(b) \leq 2(g-1)$.
  On the other hand, if $l+1 \leq g-1-|r|$
  then obviously $|2r| +2\deg(b) \leq 2(g-1)$
  as $\deg(b)\leq l+1$ by hypothesis.
\hfill $\Box$

\vspace{4mm} {\em Acknowledgments:\/} First author would like to
express his gratitude to the Fakult\"at f\"ur Mathematik and to
Prof.\ Stefan Bauer for their hospitality and support during his
stay at Universit\"at Bielefeld when part of this work was carried
out. He is specially indebted to Rogier Brussee with whom he
discussed many of the ideas that gave rise to this work. Also
thanks to Zolt\'an Szab\'o and Cliff Taubes for helpful
correspondence. Second author's research was supported by the
Australian Research Council Fellowship.

\noindent {\sc Departamento de Matem\'aticas, Facultad de
Ciencias, Universidad Aut\'onoma de Madrid, Ciudad Universitaria
Cantoblanco, 28049 Madrid, Spain} \\ {\tt \bf
vicente.munoz@uam.es}

\noindent {\sc Max-Planck-Institut f\"ur Mathematik, Vivatsgasse
7, D-53111 Bonn, Germany} \\ {\tt \bf bwang@mpim-bonn.mpg.de}

\end{document}